\documentclass[a4paper,12pt]{article}
\usepackage[a4paper,top=4.0cm,bottom=2.0cm,left=2.0cm,right=2.0cm]{geometry}

\usepackage{natbib}
\bibpunct{(}{)}{;}{d}{,}{,}

\usepackage{amsthm,amsmath,amsfonts,mathrsfs}
\usepackage{epsfig,graphicx}

\title{A Gaussian mixture ensemble transform filter\thanks{Universit\"at Potsdam, 
Institut f\"ur Mathematik, Am Neuen Palais 10, D-14469 Potsdam, Germany}}
\author{Sebastian Reich}

\begin{document}

\maketitle


\begin{abstract}
We generalize the popular ensemble Kalman filter to an ensemble transform filter where
the prior distribution can take the form of a Gaussian mixture or a
Gaussian kernel density estimator.  The design of the filter is 
based on a continuous formulation of the Bayesian filter analysis
step. We call the new filter algorithm the ensemble Gaussian mixture
filter (EGMF). The EGMF is implemented for three simple test problems 
(Brownian dynamics in one dimension, Langevin dynamics in two dimensions, and the 
three dimensional Lorenz-63 model). It is demonstrated that the EGMF
is capable to track systems with non-Gaussian uni- and multimodal ensemble distributions.
\end{abstract}


\section{Introduction}

We consider dynamical  models given in the form of ordinary differential 
equations (ODEs)
\begin{equation} \label{chap12_ode}
\dot{\bf x} = f({\bf x},t)
\end{equation}
with state variable ${\bf x} \in \mathbb{R}^N$. Initial conditions 
at time $t_0$ are not precisely known and are treated as a random variable 
instead, i.e., we assume that
\begin{equation*} 
{\bf x}(t_0) \sim \pi_0,
\end{equation*}
where $\pi_0({\bf x})$ denotes a given probability density function (PDF). 
The solution of (\ref{chap12_ode}) at time $t$ with initial condition 
${\bf x}_0$ at $t_0$ is denoted by ${\bf x}(t;t_0,{\bf x}_0)$.

The evolution of the initial PDF $\pi_0$ under the ODE (\ref{chap12_ode})
up to a time $t>t_0$ is provided by the continuity equation
\begin{equation} \label{chap12_Liouville}
\frac{\partial \pi}{\partial t} = - \nabla_{\bf x} \cdot (\pi f),
\end{equation}
which is also called Liouville's equation in the statistical mechanics 
literature \citep{sr:gardiner}. Let us denote the solution of Liouville's 
equation at observation time $t$ by $\pi({\bf x},t)$. 
In other words, solutions ${\bf x}(t;t_0,{\bf x}_0)$ with ${\bf x}_0 
\sim \pi_0$ constitute a random variable with PDF $\pi(\cdot,t)$.

For a chaotic ODE (\ref{chap12_ode}), i.e.~for an ODE with positive Lyapunov
exponents, the PDF $\pi(\cdot,t)$ will be spread out over the whole
chaotic attractor for $t \to \infty$. This in turn implies a limited solution 
predictability in the sense that the time-evolved PDF will become 
increasingly independent of the initial PDF $\pi_0$. Furthermore, even if
the initial PDF is nearly Gaussian with mean $\overline{\bf x}_0$ and small covariance matrix 
${\bf P}$, the solution ${\bf x}(t;t_0,\overline{\bf x}_0)$ will become increasingly unrepresentative
of the expectation value of the underlying random variable it is supposed to represent.

To counteract the divergence of nearby trajectories under chaotic dynamics,
we assume that we have uncorrelated measurements
${\bf y}_{\rm obs}(t_j) \in \mathbb{R}^K$ at times $t_j$, $j\ge 1$ 
with  measurement error covariance matrix 
${\bf R} \in \mathbb{R}^{K\times K}$, i.e.
\begin{equation} \label{measurement}
{\bf y}_{\rm obs}(t_j) - {\bf H}{\bf x}(t_j) \sim {\rm N}({\bf 0},{\bf R}),
\end{equation}
where the notation ${\rm N}(\overline{\bf y},{\bf B})$ is used to denote a
normal distribution in ${\bf y} \in \mathbb{R}^K$ with mean $\overline{\bf y}$ 
and covariance matrix ${\bf B} \in \mathbb{R}^{K\times K}$.
The matrix ${\bf H} \in \mathbb{R}^{K \times N}$ is called
the forward operator. The task of combining solutions to (\ref{chap12_ode}) with 
intermittent measurements (\ref{measurement}) is called data assimilation in the
geophysical literature \citep{sr:evensen}
and filtering in the statistical literature \citep{sr:crisan}.

A first step to perform data assimilation for nonlinear ODEs
(\ref{chap12_ode}) is to approximate solutions to the associated
Liouville equation (\ref{chap12_Liouville}). In this paper, we rely 
exclusively on particle methods \citep{sr:crisan} for which Liouville's 
equation is naturally approximated by the evolving empirical measure. 
More precisely, particle or ensemble filters rely on the 
simultaneous propagation of $M$ independent solutions ${\bf x}_i(t) \in 
\mathbb{R}^N$, 
$i=1,\ldots,M$, of (\ref{chap12_ode}) \citep{sr:evensen}. We associate the 
empirical measure
\begin{equation} \label{empirical}
\pi_{\rm em}({\bf x},t) = \sum_{i=1}^M \gamma_i \delta({\bf x}-{\bf x}_i(t))
\end{equation}
with weights $\gamma_i > 0$ satisfying
\begin{equation*}
\sum_{i=1}^M \gamma_i = 1.
\end{equation*}
Here $\delta(\cdot)$ denotes the Dirac delta function. Hence our statistical 
model is given by the empirical measure (\ref{empirical}) and is 
parametrized by the particle weights $\{\gamma_i\}$ and the particle 
locations $\{{\bf x}_i\}$. In the absence of measurements, the empirical 
measure $\pi_{\rm em}$ with constant weights $\gamma_i$ is an
exact (weak) solution to Liouville's equation (\ref{chap12_Liouville}) 
provided the ${\bf x}_i(t)$'s are solutions to the ODE
(\ref{chap12_ode}). Optimal statistical efficiency is achieved with
equal particle weights $\gamma_i = 1/M$.

The assimilation of a measurement at $t_j$ leads via Bayes' theorem to a discontinuous change
in the statistical model (\ref{empirical}). Sequential Monte Carlo methods \citep{sr:crisan} 
are primarily based on a discontinuous change in the weight factors $\gamma_i$. 
To avoid a subsequent degeneracy in the particle weights one
re-samples or uses other techniques which essentially lead to a redistribution of particle 
positions ${\bf x}_i$. See, for example, \cite{sr:crisan} for more
details. The ensemble Kalman filter 
(EnKF) relies on the alternative idea to replace the empirical measure (\ref{empirical}) 
by a Gaussian PDF prior to an assimilation step \citep{sr:evensen}.
This approach allows for the application of the Kalman analysis
formulas to the ensemble mean and covariance matrix. The final step
of an EnKF is the re-interpretation of the Kalman analysis step in
terms of modified particle positions while the weights are held constant
at $\gamma_i = 1/M$. We call filter algorithms that rely on modified
particle/ensemble positions and fixed particle weights {\it ensemble
transform filters}. A new ensemble transform
filter has recently been proposed by \cite{sr:anderson10}. The filter
is based on an appropriate transformation step in observation space
and subsequent linear regression of the transformation onto the full
state space. The approach developed in this paper relies instead on a general
methodology for deriving ensemble transform filters as proposed by
\cite{sr:reich10}. See Section \ref{sec_transform} below for a
summary. The same methodology has been developed for
continuous-in-time observations by \cite{sr:crisan10}. In this paper,
we demonstrate how our ensemble transform filter framework can be used
to generalize EnKFs to Gaussian mixture models and Gaussian kernel
density estimators. The essential steps
are summarized in Section \ref{sec_GM} while an algorithmic summary of
the proposed ensemble Gaussian mixture filter (EGMF) is provided in
Section \ref{sec_summary}. The EGMF can also be viewed as a generalization of the
continuous formulation of ensemble square root filters \citep{sr:tippett03}
as provided by \cite{sr:br10,sr:br10b} and the EnKF with perturbed observations, as
demonstrated by \cite{sr:reich10}. The paper concludes with three numerical
examples in Section \ref{sec_numerics}. We first demonstrate the
properties of the newly proposed EGMF for one-dimensional Brownian
dynamics under a double-well potential. This simulation is extended to the
associated two-dimensional Langevin dynamics model with only velocities being 
observed. Finally we consider the three variable model of \cite{sr:lorenz63}.

We mention that alternative extensions of the EnKF to Gaussian
mixtures have recently been proposed, for example, by \cite{sr:smith07},
\cite{sr:stordal11}, and \cite{sr:frei11}. However, while
the cluster EnKF  of \cite{sr:smith07} is an example of an ensemble 
transform filter, it fits the posterior (analysed) ensemble
distribution back to a single Gaussian PDF and, hence, only partially
works with a Gaussian mixture. Both the mixture ensemble Kalman filter of \cite{sr:frei11}
and the adaptive Gaussian mixture filter of
\cite{sr:stordal11} approximate the model uncertainty by a sum of
Gaussian kernels and utilize the ensemble Kalman filter as a particle update step
under a single Gaussian kernel. Resampling or a re-weighting of particles is required to avoid a degeneracy 
of particle weights due to changing kernel weights. A related filter algorithm based 
on Gaussian kernel density estimators has previously been considered by \cite{sr:andand99}.

\section{A general framework for ensemble transform filters} \label{sec_transform}

Bayes' formula can be interpreted as a discontinuous change of a forecast PDF $\pi_f$ 
into an analyzed PDF $\pi_a$ at each observation time $t_j$. On the other hand, one can
find a continuous embedding $\pi({\bf x},s)$ with respect to 
a fictitious time $s \in [0,1]$ such that $\pi(\cdot,0) = \pi_f$ and
$\pi_a = \pi(\cdot,1)$. As proposed by \cite{sr:reich10}, 
this embedding can be viewed as being induced by a continuity (Liouville) equation
\begin{equation} \label{optimaltransportation1}
\frac{\partial \pi}{\partial s} = -\nabla_{\bf x} \cdot \left(\pi g  \right)
\end{equation}
for an appropriate vector field $g({\bf x},s) \in \mathbb{R}^N$. The vector field $g$ is not
uniquely determined for a given continuous embedding $\pi(\cdot,s)$ unless we also require that 
it is the minimizer of the kinetic energy
\begin{equation*}
{\cal T}(v) = \frac{1}{2} \int {\rm d}\pi  \, v^T {\bf M} v
\end{equation*}
over all admissible vector fields $v \in L^2({\rm d}\pi,\mathbb{R}^N)$, where ${\bf M} \in 
\mathbb{R}^{N\times N}$ is a positive definite mass matrix \citep{sr:Villani}. 
Admissibility means that $g=v$ satisfies
(\ref{optimaltransportation1}) for given $\pi$ and $\partial \pi/\partial s$.

Under these assumptions,  a constrained variational principle \citep{sr:Villani} implies that 
the desired vector field is given by $g = {\bf M}^{-1} \nabla_{\bf x} \psi$, where 
the potential $\psi({\bf x},s)$ is the solution of the elliptic partial differential equation (PDE)
\begin{equation} \label{optimaltransportation2}
\nabla_{\bf x} \cdot \left( \pi {\bf M}^{-1} \nabla_{\bf x} \psi \right) = 
\pi \left( S - \mathbb{E}_\pi[S]\right)
\end{equation}
for given PDF $\pi$, mass matrix ${\bf M}$, and negative log-likelihood function 
\begin{equation} \label{nlogl}
S({\bf x};{\bf y}_{\rm obs}(t_j)) = \frac{1}{2}
({\bf H}{\bf x}-{\bf y}_{\rm obs}(t_j))^T {\bf R}^{-1} ({\bf H}{\bf x}-{\bf y}_{\rm obs}
(t_j)).
\end{equation}
Here $\mathbb{E}_\pi[f]$ denotes the expectation value of a function $f({\bf x})$ with
respect to a PDF $\pi({\bf x})$. We finally replace (\ref{optimaltransportation1}) by
\begin{equation} \label{alphaODEIII}
\frac{\partial \pi}{\partial s}  = -\nabla_{\bf x} \cdot \left( \pi 
{\bf M}^{-1} \nabla_{\bf x} \psi \right)
\end{equation}
with an underlying ODE formulation
\begin{equation} \label{alphaODEIV}
\frac{{\rm d} {\bf x}}{{\rm d} s} = {\bf M}^{-1} \nabla_{\bf x} \psi({\bf x},s)
\end{equation}
in fictitious time $s \in [0,1]$. As for the ODE (\ref{chap12_ode}) and its associated Liouville
equation (\ref{chap12_Liouville}), we may approximate (\ref{alphaODEIV}) and its associated 
Liouville equation (\ref{alphaODEIII}) by an empirical measure of type (\ref{empirical}).
Furthermore, one and the same empirical measure approximation can now be used for both
the ensemble propagation step under the model dynamics (\ref{chap12_ode}) and the
data assimilation step (\ref{alphaODEIII}) using constant and equal weights $\gamma_i = 1/M$.
The particle filter approximation is closed by finding an appropriate numerical solution to
the elliptic PDE (\ref{optimaltransportation2}). This is the crucial step which will lead to
different ensemble transform filter algorithms. 

The basic numerical approach to the data assimilation step within an ensemble transform filter 
formulation consists then of the following sequence 
of steps. (i) Given a current ensemble of solutions ${\bf x}_i(s)$,
$i=1,\ldots,M$, one fits a statistical model $\widehat{\pi}({\bf x},s)$. (ii) Solve the elliptic PDE
\begin{equation} \label{spde}
\nabla_{\bf x} \cdot (\widehat{\pi} \widehat{g} ) = \widehat{\pi} \left(S - \mathbb{E}_{\widehat{\pi}}[S] \right)
\end{equation}
for a vector field $\widehat{g}({\bf x},s)$. The solution is not
uniquely determined and an appropriate choice needs to be made. See
the discussion above. (iii) Propagate the ensemble members under the ODE
\begin{equation} \label{sode}
\frac{{\rm d}  {\bf x}_i}{{\rm d}s} = \widehat{g}({\bf x}_i,s).
\end{equation}
We assume that a forecast ensemble of
$M$ members ${\bf x}_i \in \mathbb{R}^N$, $i=1,\ldots,M$, is available at
an observation time $t_j$ which provides the initial conditions for the ODE (\ref{sode}). 
Solutions at $s=1$ yield the analyzed ensemble members, 
which are then used as the new initial conditions for (\ref{chap12_ode}) at
time $t=t_j$ and (\ref{chap12_ode}) is solved over $[t_j,t_{j+1}]$ up to the next observation point.

If the statistical model is a Gaussian with mean $\overline{\bf x}\in \mathbb{R}^N$ and covariance matrix
${\bf P} \in \mathbb{R}^{N\times N}$, then the outlined approach leads to a continuous formulation of the 
ensemble square-root ensemble filter analysis step at time $t_j$  \citep{sr:br10,sr:br10b}, i.e.
\begin{equation} \label{enkf}
\frac{{\rm d}{\bf x}_i}{{\rm d}s} =
-\frac{1}{2}  {\bf P} {\bf H}^T {\bf R} ^{-1} \left(
{\bf H} {\bf x}_i + {\bf H} \overline{\bf x} - 2{\bf y}_{\rm obs}(t_j)\right)
\end{equation}
for $s\in [0,1]$. It follows that ${\bf M} = {\bf P}^{-1}$ and
\begin{equation}
\psi({\bf x}) = -\frac{1}{4} \left({\bf H}{\bf x} + {\bf H}\overline{\bf
    x} - 2{\bf y}_{\rm obs}(t_j) \right)^T {\bf R}^{-1} 
\left({\bf H}{\bf x} + {\bf H}\overline{\bf
    x} - 2{\bf y}_{\rm obs}(t_j) \right).
\end{equation}

\section{An ensemble transform filter based on 
Gaussian mixture statistical models} \label{sec_GM}

We now develop an ensemble transform filter algorithm based on a $L\ge
1$ component Gaussian mixture model, i.e.
\begin{equation} \label{mixturem}
\widehat{\pi}({\bf x}) = \sum_{l=1}^L \frac{\alpha_{l}}{(2\pi)^{N/2}
\det {\bf P}_l^{1/2}} 
\exp \left(-\frac{1}{2} ({\bf x}-\overline{\bf x}_{l})^T {\bf
  P}_l^{-1} ({\bf x}-\overline{\bf x}_{l})
 \right) = \sum_{l=1}^L \alpha_{l} \pi_{{\rm Gauss},l}({\bf x}),
\end{equation}
where $\pi_{{\rm Gauss},l}({\bf x})$ denotes the normal distribution
in ${\bf x}\in \mathbb{R}^N$ with mean $\overline{\bf x}_l$ and 
covariance matrix ${\bf P}_l$.
The Gaussian mixture parameters, i.e.~$\alpha_{l}, \overline{\bf x}_{l},
{\bf P}_{l}$, $l=1,\ldots,L$, need to be determined from the ensemble
${\bf x}_i$, $i=1,\ldots,M$, in an appropriate manner
using, for example, the expectation-maximization (EM) algorithm
\citep{sr:dempster,sr:smith07}. See Section \ref{sec3.3} for more details.
Note that $\sum_{l=1}^L \alpha_l = 1$ and $\alpha_l \ge 0$.
To simplify notation, we write $\pi_l$ instead of $\pi_{{\rm Gauss},l}$
from now on. 

An implementation of the associated continuous formulation of the
Bayesian analysis step proceeds as follows. To simplify the discussion we derive our filter formulation
for a scalar observation variable, i.e.~$K=1$, $y_{\rm obs}(t_j) - {\bf H}{\bf x}(t_j)  \sim {\rm N}(0,R)$, and
\begin{equation}
S({\bf x};y_{\rm obs}(t_j)) = \frac{1}{2R} \left( y_{\rm obs}(t_j) - {\bf H}{\bf x}\right)^2 .
\end{equation}
The vector-valued case can be treated accordingly provided the error covariance matrix ${\bf R}$ is diagonal. 
We first decompose the vector
field $\widehat{g}({\bf x},s) \in \mathbb{R}^N$ in (\ref{sode}) into two contributions, i.e.
\begin{equation} \label{total}
\frac{{\rm d} {\bf x}}{{\rm d}s} = \widehat{g}({\bf x},s) = {u}_{\rm A}({\bf x},s) + {u}_{\rm B} ({\bf x},s) .
\end{equation}
To simplify notation we drop the explicit $s$ dependence in the following calculations.
We next decompose the right hand side of the elliptic PDE (\ref{spde}) also into two contributions
\begin{equation}
\widehat{\pi} \left(S({\bf x}) -  \mathbb{E}_{\widehat{\pi}}[S] \right) 
= \left\{\sum_{l=1}^L \alpha_l \pi_l  \left( S({\bf x}) - \mathbb{E}_{\pi_l}[S]\right) \right\}
+ \left\{ \sum_{l=1}^L \alpha_l \pi_l \left( \mathbb{E}_{\pi_l} - \mathbb{E}_{\widehat{\pi}}[S] \right) \right\}.
\end{equation}
We now derive explicit expressions for $u_{\rm A}({\bf x})$ and $u_{\rm B}({\bf u})$.

\subsection{Gaussian mixture Kalman filter contributions}

We define the vector field ${u}_{\rm A}({\bf x})$ through the equation
\begin{equation} \label{ansatzA}
{u}_{\rm A} ({\bf x}) = \sum_{l=1}^L \frac{\alpha_l
  \pi_{l}({\bf x})}{\widehat{\pi}({\bf x})} {\bf P}_l
\nabla_{\bf x} \psi_{{\rm A},l} ({\bf x}),
\end{equation}
together with
\begin{equation}
\nabla_{\bf x} \cdot \left\{ \widehat{\pi}({\bf x}) {u}_{\rm
    A}({\bf x}) \right\} =
\nabla_{\bf x} \cdot \left\{ \sum_{l=1}^L \alpha_l \pi_l({\bf x}) {\bf P}_l \nabla_{\bf x} \psi_{{\rm
      A},l}({\bf x}) \right\} = \sum_{l=1}^L \alpha_l \pi_l({\bf x}) (S({\bf x}) -
    \mathbb{E}_{\pi_l}[S] )
\end{equation}
which implies that the potentials $\psi_{{\rm A},l}({\bf x})$,
$l=1,\ldots,L$, are uniquely determined by
\begin{equation}
\nabla_{\bf x} \cdot \left\{  \pi_l({\bf x}) {\bf P}_l \nabla_{\bf x} \psi_{{\rm
      A},l}({\bf x}) \right\} =  \pi_l({\bf x}) (S({\bf x}) -
    \mathbb{E}_{\pi_l}[S] )
\end{equation}
for all $l=1,\ldots,L$. It follows that the potentials $\psi_{{\rm
    A},l}({\bf x})$ are equivalent to the ensemble Kalman filter potentials for the $l$-th
Gaussian component. Hence, using (\ref{enkf}) and (\ref{ansatzA}), we obtain
\begin{equation} \label{uA}
{u}_{\rm A}({\bf x},s) = -\frac{1}{2} 
  \sum_{l=1}^L \frac{\alpha_l(s) \pi_l({\bf x},s)}{\widehat{\pi}({\bf x},s)} {\bf P}_l(s) {\bf H}^T {R}^{-1}
\left[ {\bf H}  {\bf x}(s) + {\bf H} \overline{\bf x}_l(s) - 2y_{\rm
    obs}(t_j) \right].
\end{equation}

\subsection{Gaussian mixture exchange contributions}

The remaining contributions for solving (\ref{optimaltransportation1})
are collected in the vector field 
\begin{equation}
{u}_{\rm B} ({\bf x}) = \sum_{l=1}^L \frac{\alpha_l
  \pi_{l}({\bf x})}{\widehat{\pi}({\bf x})} 
{\bf P}_l \nabla_{\bf x} \psi_{{\rm B},l} ({\bf x}),
\end{equation}
which therefore needs to satisfy
\begin{equation}
\nabla_{\bf x} \cdot \left\{ \widehat{\pi}({\bf x}) {u}_{\rm
    B}({\bf x}) \right\} =
\nabla_{\bf x} \cdot \left\{ \sum_{l=1}^L \alpha_l \pi_l({\bf x}) {\bf P}_l \nabla_{\bf x} \psi_{{\rm
      B},l}({\bf x}) \right\} = \sum_{l=1}^L \alpha_l \pi_l({\bf x})
(\mathbb{E}_{\pi_l} [S] -
    \mathbb{E}_{\widehat{\pi}}[S] )
\end{equation}
and, hence, we may set
\begin{equation} \label{b}
\nabla_{\bf x} \cdot \left\{  \pi_l({\bf x}) {\bf P}_l \nabla_{\bf x} \psi_{{\rm
      B},l}({\bf x}) \right\} =  \pi_l({\bf x}) (\mathbb{E}_{\pi_l} [S] -
    \mathbb{E}_{\widehat{\pi}}[S] )
\end{equation}
for all $l=1,\ldots,L$. 
To find a solution of (\ref{b}) we introduce functions $\widehat{\psi}_{{\rm B},l}$ such that
\begin{equation}
\psi_{{\rm B},l}({\bf x}) = \widehat{\psi}_{{\rm B},l}({\bf H}{\bf x}-{\bf H}\overline{\bf x}_l) = 
\widehat{\psi}_{{\rm B},l}(y-\overline{y}_l)
\end{equation}
with $y := {\bf H}{\bf x}$ and $\overline{y}_l = {\bf H}\overline{\bf x}_l$. 
Now the right hand side of (\ref{b}) gives rise to
\begin{equation}
\nabla_{\bf x} \cdot \left\{  \pi_l({\bf x}) {\bf P}_l \nabla_{\bf x} \psi_{{\rm
      B},l}({\bf x}) \right\}
= \pi_l({\bf x}) \left( - (y-\overline{y}_l) \frac{{\rm d} \widehat{\psi}_{{\rm
      B},l}}{{\rm d}y} (y-\overline{y}_l) + {\bf H}{\bf P}_l {\bf H}^T \frac{{\rm d}^2 \widehat{\psi}_{{\rm
      B},l}}{{\rm d}y^2} (y-\overline{y}_l) \right)
\end{equation}
and (\ref{b}) simplifies further to
the scalar PDE
\begin{equation} \label{bb}
-(y-\overline{y}_l) \frac{{\rm d} \widehat{\psi}_{{\rm
      B},l}}{{\rm d}y} (y-\overline{y}_l) + {\bf H}{\bf P}_l {\bf H}^T
\frac{{\rm d}^2 \widehat{\psi}_{{\rm
      B},l}}{{\rm d}y^2} (y-\overline{y}_l)
 = \mathbb{E}_{\pi_l} [S] - \mathbb{E}_{\widehat{\pi}}[S] .
\end{equation}
The PDE (\ref{bb}) can be solved for 
\begin{equation}
f(z) = \frac{{\rm d} \widehat{\psi}_{{\rm
      B},l}}{{\rm d}y} (y-\overline{y}_l), \qquad z = y - \overline{y}_l,
\end{equation}
under the condition $f(0) = 0$ by explicit quadrature and one obtains
\begin{equation}
f(y-\overline{y}_l) = \frac{1}{2} \frac{\mathbb{E}_{\pi_l} [S] - \mathbb{E}_{\widehat{\pi}}[S] 
}{{\bf H} {\bf P}_l {\bf H}}  \frac{ \mbox{erf} \left(  (y-\overline{y}_l)/ \sqrt{2 \sigma_l^2 }
\right) }{\pi_l(y)}
\end{equation}
with marginalized PDF
\begin{equation} \label{md}
\pi_l(y) := \frac{1}{ \sqrt{2\pi  \sigma_l^2} } \exp \left(
-\frac{\left( y - \overline{y}_l\right)^2}{2 \sigma_l^2} \right) ,
\end{equation}
$\sigma_l = \sqrt{{\bf H}{\bf P}_l {\bf H}^T}$,
and the standard error function
\begin{equation}
{\rm erf}(y) =  \frac{2}{\sqrt{\pi}} \int_0^y e^{-s^2}{\rm d}s .
\end{equation}
Note that
\begin{equation} \label{deriv}
\frac{1}{2}  \frac{{\rm d}}{{\rm d}y} \mbox{erf}\left((y-\overline{y}_l)/ \sqrt{2 \sigma_l^2}
\right) = \pi_l(y).
\end{equation}
We finally obtain the expression
\begin{equation} \label{uB}
{u}_{\rm B} ({\bf x},s) = \frac{1}{2} \sum_{l=1}^L \frac{\alpha_l(s)
  \pi_{l}({\bf x},s)}{\widehat{\pi}({\bf x},s)}   {\bf P}_l(s) {\bf
  H}^T \frac{ \mathbb{E}_{\pi_l}[S](s)- \mathbb{E}_{\widehat{\pi}}[S](s)  }{\sigma_l^2}
\frac{\mbox{erf}\left((y-\overline{y}_l)/\sqrt{2 \sigma_l^2}
\right) }{\pi_l(y)}
\end{equation}
for the vector field $u_{\rm B}({\bf x},s)$. 

The expectation values $\mathbb{E}_{\pi_{l}}[S]$, $l=1,\ldots,L$, 
can be computed analytically, i.e.
\begin{equation}
\mathbb{E}_{\pi_{l}}[S] = \frac{1}{2R}\left( (y_{\rm obs}(t_j) - \overline{y}_l)^2 + \sigma_l^2 \right)
\end{equation}
or estimated numerically. 
Recall that $\sum \alpha_l  = 1$ and, therefore,
\begin{equation}
\mathbb{E}_{\widehat{\pi}}[S] = \sum_{l=1}^L \alpha_l \mathbb{E}_{\pi_{l}}[S] .
\end{equation}

It should be kept in mind that the Gaussian mixture parameters
$\alpha_l, \overline{{\bf x}}_l,{\bf P}_l$ can be updated directly using the
differential equations
\begin{eqnarray}
\frac{{\rm d} \overline{\bf x}_l}{{\rm d}s} &=& -
{\bf P}_l   {\bf H}^T {R}^{-1} ({\bf H} \overline{{\bf x}}_l -
{\bf y}_{\rm obs}(t_j)),\\
\frac{{\rm d} {\bf P}_l}{{\rm d}s} &=& - {\bf P}_l  {\bf H}^T {R}^{-1}
{\bf H}  {\bf P}_l,\\
\frac{{\rm d} \alpha_l}{{\rm d}s} &=& - \frac{1}{2}\alpha_l \left\{ 
({\bf H} \overline{\bf x}_l- {\bf y}_{\rm obs}(t_j))^T {R}^{-1} 
({\bf H} \overline{\bf x}_l- {\bf y}_{\rm obs}(t_j)) + \lambda \right\}, \label{alpha}
\end{eqnarray}
for $l=1,\ldots,L$. Here $\lambda \in \mathbb{R}$ is chosen such that
\begin{equation}
\sum_{l=1}^L \frac{{\rm d}\alpha_l}{{\rm d} s} = 0.
\end{equation}
Furthermore, ${u}_{\rm A}({\bf x},s)$ exactly mirrors the update
of the Gaussian mixture means $\overline{\bf x}_l$ and covariance
matrices ${\bf P}_l$, while ${u}_{\rm  B}({\bf x},s)$ mimics the
update of the weight factors $\alpha_l$ by rearranging the particle positions
accordingly. 

As already eluded to, we can treat each
uncorrelated observation separately and sum the individual contributions in $u_{\rm
  A}({\bf x},s)$ and $u_{\rm B}({\bf x},s)$, respectively, to obtain the desired total vector
field (\ref{total}). 

\subsection{Implementation aspects} \label{sec3.3}

Given a set of ensemble members ${\bf x}_i$, $i=1,\ldots,M$, there are several options for implementing
a Gaussian mixture filter. The trivial case $L=1$ leads back to the
continuous formulations of \cite{sr:br10,sr:br10b}. More
interestingly, one can chose $L\ll M$ and estimate the mean and the covariance 
matrices for the Gaussian mixture model using the EM algorithm \citep{sr:dempster,sr:smith07}. 
The EM algorithm self-consistently computes the mixture mean $\overline{\bf x}_l$ and covariance matrix 
${\bf P}_l$ via
\begin{equation} \label{emmc}
\overline{\bf x}_l = \frac{1}{\sum_{i=1}^M \beta_{i,l}} \sum_{i=1}^M \beta_{i,l} {\bf
    x}_i, \qquad
{\bf P}_l = \frac{1}{\sum_{i=1}^M \beta_{i,l}} \sum_{i=1}^M \beta_{i,l}
  \left( {\bf x}_i - \overline{\bf x}_l\right) \left( {\bf x}_i -
    \overline{\bf x}_l \right)^T 
\end{equation}
for $l=1,\ldots,L$ using weights $\beta_{i,l}$ defined by 
\begin{equation} \label{EM1}
\beta_{i,l} = \frac{\alpha_l \pi_l ({\bf x}_i)}{\sum_{k=1}^M 
\alpha_k \pi_k({\bf x}_i)} , \qquad \alpha_l = \frac{1}{M} \sum_{i=1}^M \beta_{i,l} .
\end{equation}
The EM algorithm can
fail to converge and a possible remedy is to add a constant contribution $\delta
{\bf I}$  to the empirical covariance matrix in (\ref{emmc}) with the
parameter $\delta >0$ appropriately chosen. We mention that more refined
implementations of the EM algorithm, such as those discussed by
\cite{sr:fraley07}, could also be considered. It is also possible to select the number of mixture
components adaptively. See, for example, \cite{sr:smith07}. 
The resulting vector fields for the $i$th ensemble member are given by
\begin{equation} \label{EGM1}
{u}_{\rm A}({\bf x}_i,s) = -\frac{1}{2} \sum_{l=1}^L \beta_{i,l}(s) {\bf P}_l(s) {\bf H}^T {R}^{-1}
\left[ {\bf H}  {\bf x}_i(s)  + {\bf H} \overline{\bf x}_l (s) - 2y_{\rm
    obs}(t_j)  \right]
\end{equation}
and, using (\ref{uB}),
\begin{equation} \label{EGM2}
u_{\rm B}({\bf x}_i,s) = \frac{1}{2}
\sum_{l=1}^L \beta_{i,l}(s)   {\bf P}_l(s) {\bf
  H}^T \frac{ \mathbb{E}_{\pi_l}[S](s)- \mathbb{E}_{\widehat{\pi}}[S](s)  }{\sigma_l^2}
\frac{\mbox{erf}\left((y-\overline{y}_l)/\sqrt{2 \sigma_l^2}
\right) }{\pi_l(y)}
\end{equation}
with weights $\beta_{i,l}$ given by (\ref{EM1}).

Another option to implement an EGMF is to set the number of mixture components equal to the number of
ensemble members, i.e.~$L=M$, and to use a prescribed covariance matrix ${\bf B}$ for all mixture
components, i.e.~${\bf P}_l = {\bf B}$ and $\overline{\bf x}_l = {\bf x}_l$, $l=1,\ldots,L$. We also give all
mixture components equal weights $\alpha_l = 1/M$. In this setting, it is more appropriate to 
call (\ref{mixturem}) a kernel density estimator \citep{sr:wand}. Then 
\begin{equation} \label{MA}
{u}_{\rm A}({\bf x}_i,s) = -\frac{1}{2} 
  \sum_{l=1}^L \beta_{i,l}(s) {\bf B} {\bf H}^T {R}^{-1}
\left[ {\bf H}  {\bf x}_i(s) + {\bf H} {\bf x}_l(s) - 2y_{\rm
    obs}(t_j) \right]
\end{equation}
and
\begin{equation} \label{FGM1}
u_{\rm B}({\bf x}_i,s) = \frac{1}{2}
\sum_{l=1}^L \beta_{i,l}(s)   {\bf B} {\bf
  H}^T \frac{ \mathbb{E}_{\pi_l}[S](s)-
  \mathbb{E}_{\widehat{\pi}}[S](s)  }{{\bf H} {\bf B} {\bf H}^T}
\frac{\mbox{err}\left((y-\overline{y}_l)/\sqrt{2 \sigma_l^2}
\right) }{\pi_l(y)}
\end{equation}
with weights $\beta_{i,l}$ given by (\ref{EM1}), $\sigma_l = \sqrt{{\bf H}{\bf B} {\bf H}^T}$, and $\alpha_l  = 1/M$. 
The Kalman filter like contributions (\ref{MA}) can be replaced by a formulation with
perturbed observations \citep{sr:evensen,sr:reich10} which yields
\begin{equation} \label{FGM2}
{u}_{\rm A}({\bf x}_i,s) = -{\bf B} {\bf H}^T {R}^{-1}
\left[ {\bf H}  {\bf x}_i(s)  - y_{\rm
    obs}(t_j) + d_i \right],
\end{equation}
where $d_i \in \mathbb{R}$, $i=1,\ldots,m$, are independent, 
identically distributed Gaussian random numbers
with mean zero and variance $R$.  A particular choice is ${\bf B} = c{\bf P}$, where ${\bf P}$ is the empirical 
covariance matrix of the ensemble and $c>0$ is an appropriate
constant. Assuming that the underlying  probability
density is Gaussian with covariance ${\bf P}$, the choice
\begin{equation} \label{bandwidth}
c = (2/(N+2))^{4/(N+4)} M^{-2/(N+4)}
\end{equation}
is optimal for large ensemble sizes $M$ in the sense of kernel density estimation (see, e.g.~\cite{sr:wand}). 
Recall that $N$ denotes the dimension of phase space. 
The resulting filter is then similar in spirit to the rank histogram filter (RHF)
suggested by \cite{sr:anderson10} with the RHF increments in observation space being replaced by those
generated from a Gaussian kernel density estimator. Another choice is ${\bf B} \approx {\bf P}$ in which
case (\ref{FGM1}) can be viewed as a correction term to the standard ensemble Kalman filter (\ref{FGM2}). 
We will explore the kernel estimator in the numerical experiment of Section \ref{sec_num3}. 

Note that localization, as
introduced by \cite{sr:houtekamer01} and \cite{sr:hamill01}, can be
combined with (\ref{EGM1})-(\ref{EGM2}) and (\ref{MA})-(\ref{FGM1}), respectively, 
as outlined in \cite{sr:br10}.  For example, one could set the covariance matrix ${\bf B}$ in
(\ref{MA})-(\ref{FGM1}) equal to the localized ensemble covariance matrix. 
Localization will be essential for a successful application of the
proposed filter formulations to high-dimensional systems (\ref{chap12_ode}).
The same applies to ensemble inflation \citep{sr:andand99}.

We also note that the computation of the particle-mixture weight factors
(\ref{EM1}) can be become expensive since it requires the computation of ${\bf
  P}_l^{-1}$. This can be avoided by either using only the diagonal part of
${\bf P}_l$ in $\pi_l({\bf x}_i)$ \citep{sr:smith07} or  by using a
marginalized density such as (\ref{md}), i.e.~$\pi_l(y_i)$, $y_i := {\bf H}{\bf x}_i$, instead
of the full Gaussian PDF values $\pi_l({\bf x}_i)$. Some other suitable
marginalization could also be performed.  

The vector field $u_{\rm B}({\bf x},s)$ is responsible for a transfer of
particles between different mixture components according to the
observation adjusted relative weight $\alpha_l$ of each mixture component. 
These transitions can be rather rapid implying that $\|u_{\rm B}({\bf
  x},s)\|_\infty$ can become large in magnitude. This can pose numerical
difficulties and we eliminated those by limiting the
$l_\infty$-norm of $u_{\rm B}({\bf  x},s)$ through a cut-off
value $u_{\rm cut}$. Alternatively, we might want to replace 
(\ref{md}) by a PDF which leads to less stiff contributions to
the vector field $u_{\rm B}({\bf x},s)$ such as the Student's t-distributions \citep{sr:Schaefer}.
Hence a natural approximative PDF is provided by the the scaled t-distribution with three degrees of 
freedom, i.e.
\begin{equation}
\phi(y;\bar y,\sigma) = \frac{2\sigma^3}{\pi} \frac{1}{\left( \sigma^2 + (y-\bar y)^2 \right)^2}.
\end{equation}
We also introduce the shorthand $\phi_l(y) = \phi(y;\bar
y_l,\sigma_l)$ with $\bar y_l = {\bf H}\overline{\bf x}_l$ and
$\sigma_l = \sqrt{{\bf H}{\bf P}_l{\bf H}^T}$. A first observation is that
\begin{equation}
\mathbb{E}_{\phi_l}[y] = \bar y_l , \qquad \mathbb{E}_{\phi_l}[(y-\bar y_l)^2] = \sigma_l^2 ,
\end{equation}
i.e., the first two moments of $\phi_l$ match those of (\ref{md}). The second observation is that
$\phi_l$ can be integrated explicitly, i.e.
\begin{equation}
\Phi_l(y) = \int_{\bar y_l}^y \phi_l(u) {\rm d}u = \frac{1}{\pi}
\arctan \left( \frac{(y-\bar y_l)}{\sigma_l} \right) +
\frac{\sigma_l}{\pi} \frac{(y-\bar y_l)}{\sigma_l^2 + (y-\bar y_l)^2}.
\end{equation}
Hence the relation (\ref{deriv}) suggests the alternative formulation
\begin{equation} \label{variantB}
u_{\rm B}({\bf x},s) = 
\sum_{l=1}^L \frac{\alpha_l(s)
  \pi_{l}({\bf x},s)}{\widehat{\pi}({\bf x},s)}   {\bf P}_l(s) {\bf
  H}^T \frac{ \mathbb{E}_{\pi_l}[S](s)- \mathbb{E}_{\widehat{\pi}}[S](s)  }{\sigma_l^2}
\frac{\Phi_l(y,s)) }{\phi_l(y,s)} ,
\end{equation}
where $\phi_l(y,s) = \phi (y;\bar y_l(s),\sigma_l(s))$. 

The differential equation (\ref{total}) needs to be approximated by a
numerical time-stepping procedure. In this paper, we use the forward
Euler method for simplicity. However, the limited region of stability
of an explicit method such as forward Euler
implies that the step-size $\Delta s$ needs to be made sufficiently
small. This issue has been investigated by \cite{sr:amezcua} for the 
formulation with $L=1$ (standard continuous ensemble Kalman filter
formulation) and a diagonally implicit scheme has been proposed which
overcomes the stability restrictions of the forward Euler method while 
introducing negligible computational overhead. The computational cost
of a single evaluation of (\ref{total}) for given mixture components 
should be slightly lower than for a single EnKF step since no matrix
inversion is required. Additional expenses arise from fitting the
mixture components (e.g.~using the EM algorithm) and from having to
use a number of time-steps $1/\Delta s > 1$.

\section{Algorithmic summary of the ensemble Gaussian mixture filter (EGMF)}
\label{sec_summary}

Since the proposed methodology for treating nonlinear filter problems
is based on an extension of the EnKF approach, we call the new filter the
ensemble Gaussian mixture filter (EGMF). We now provide an algorithmic summary.

First a set of $M$ ensemble members ${\bf x}_i(0)$ is generated at time
$t_0$ according to the initial PDF $\pi_0$. 

In between observations, the ensemble members are propagated under
the ODE model (\ref{chap12_ode}), i.e.
\begin{equation} \label{ensembleODE}
\dot{\bf x}_i = f({\bf x}_i,t),
\end{equation}
for $i = 1,\ldots,M$ and $t \in [t_{j-1},t_{j}]$. 

At an observation time $t_j$, a Gaussian mixture model (\ref{mixturem}) 
is fitted to the ensemble members
${\bf x}_i$, $i = 1,\ldots,M$. One can, for example, use the classic
EM algorithm \citep{sr:dempster,sr:smith07} for this purpose. In this paper we use a simple
heuristic to determine the  number of components $L \in \{1,2\}$. An adaptive selection of $L$ is, 
however, feasible (see, e.g., \cite{sr:smith07}). Alternatively,
one can set $L=M$ and implement the EGMF with a
Gaussian kernel density estimator with $\overline{\bf x}_l = {\bf
  x}_l$, $\alpha_l = 1/M$. The covariance matrix ${\bf B}$ can be 
based on the empirical covariance matrix ${\bf P}$ of the whole
ensemble via ${\bf B} = c {\bf P}$ with the constant $c>0$
appropriately chosen. At this stage covariance localization can also be applied. 

The vector fields $u_{\rm A}({\bf x},s)$ and $u_{\rm B}({\bf x},s)$
are computed according to (\ref{EGM1}) and (\ref{EGM2}), respectively, 
(or, alternatively, use (\ref{FGM1})-(\ref{FGM2})) for each independent observation 
and the resulting contributions are summed up to a total vector field $\widehat{g}({\bf x},s)$.
Next the ensemble members are updated according to
(\ref{total}) for ${\bf x} = {\bf x}_i$, $i=1,\ldots,M$. 
Here we use a simple forward Euler discretization with
step-size $\Delta s$ chosen sufficiently small. 
After each time-step the Gaussian mixture components are updated, if
necessary, using the EM algorithm.  The analyzed ensemble
members ${\bf x}_i(t_j)$ are obtained after $1/\Delta s$ time-steps as
the numerical solutions at $s=1$ and provide the new initial
conditions for (\ref{ensembleODE}) with time $t$ now in the interval $[t_j,t_{j+1}]$. 

Ensemble induced estimates for the expectation value of a function 
$f({\bf x})$ can be computed via
\begin{equation}
\overline{f} = \frac{1}{M} \sum_{i=1}^M f({\bf x}_i)
\end{equation}
without reference to a Gaussian mixture model.

\section{Numerical experiments} \label{sec_numerics}

In this section we provide results from several numerical simulations 
and demonstrate the performance of the proposed EGMF in comparison
with standard implementations of the EnKF and an implementation
of the RHF \citep{sr:anderson10}. We first  investigate the
Bayesian assimilation step without model equations.

\subsection{Single Bayesian assimilation step}

\begin{figure}
\begin{center}
\includegraphics[width=0.4\textwidth]{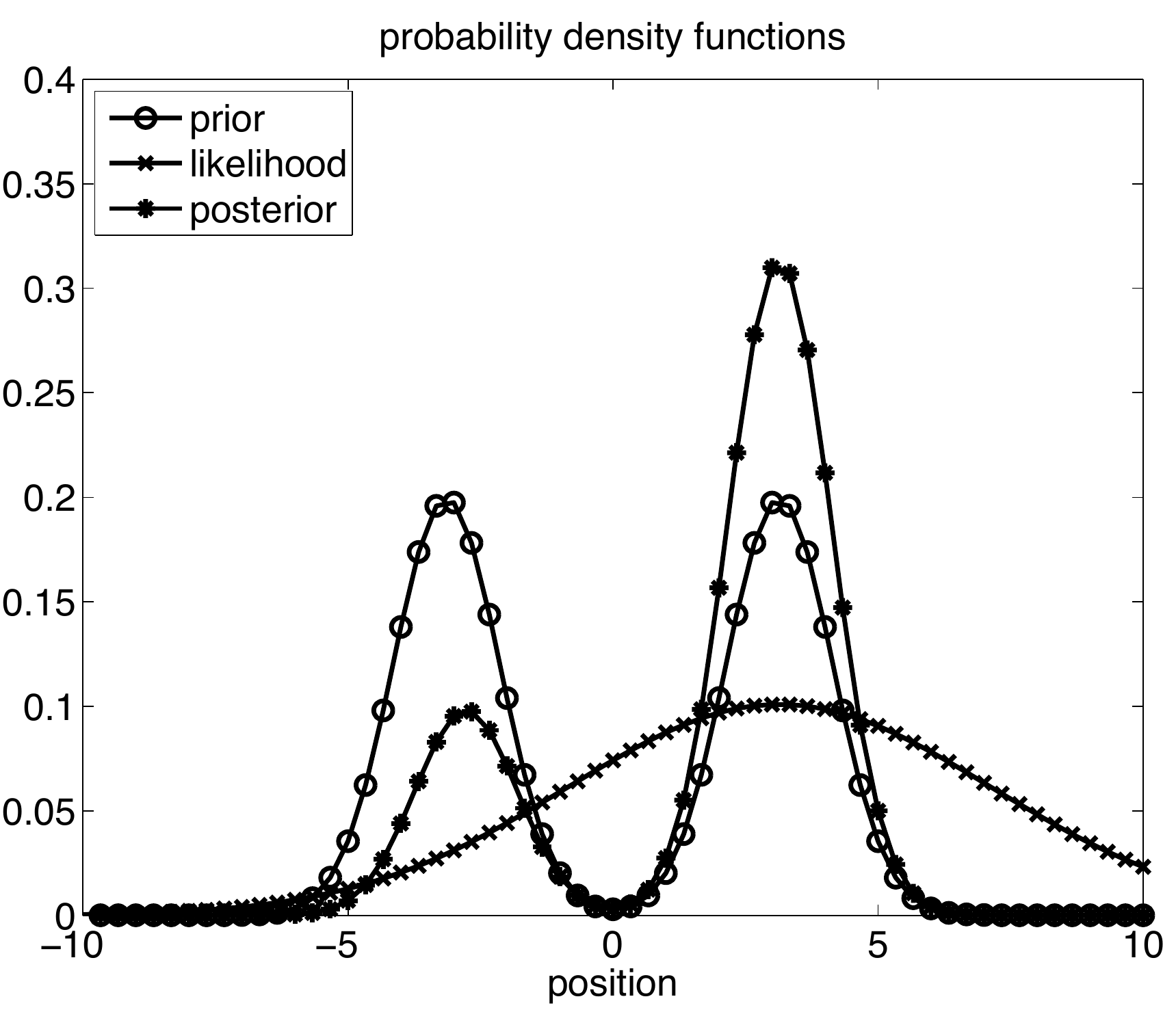}
\end{center}
\caption{Displayed are the prior distribution, the likelihood from a
  measurement and the resulting posterior distribution. The prior as
  well as the posterior are bimodal Gaussian.}
\label{figN1.1}
\end{figure}

\begin{figure}
\begin{center}
\includegraphics[width=0.4\textwidth]{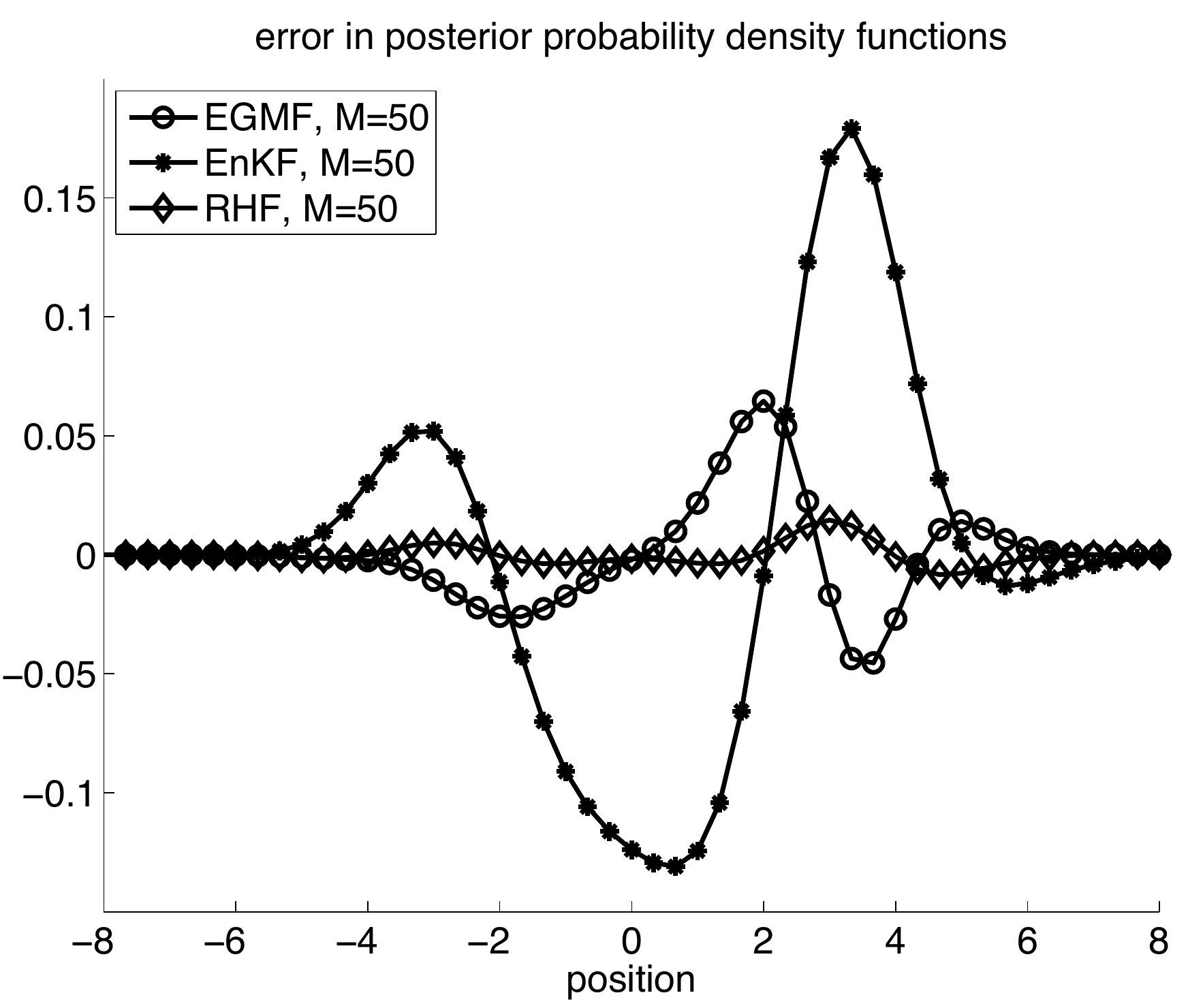}
\includegraphics[width=0.4\textwidth]{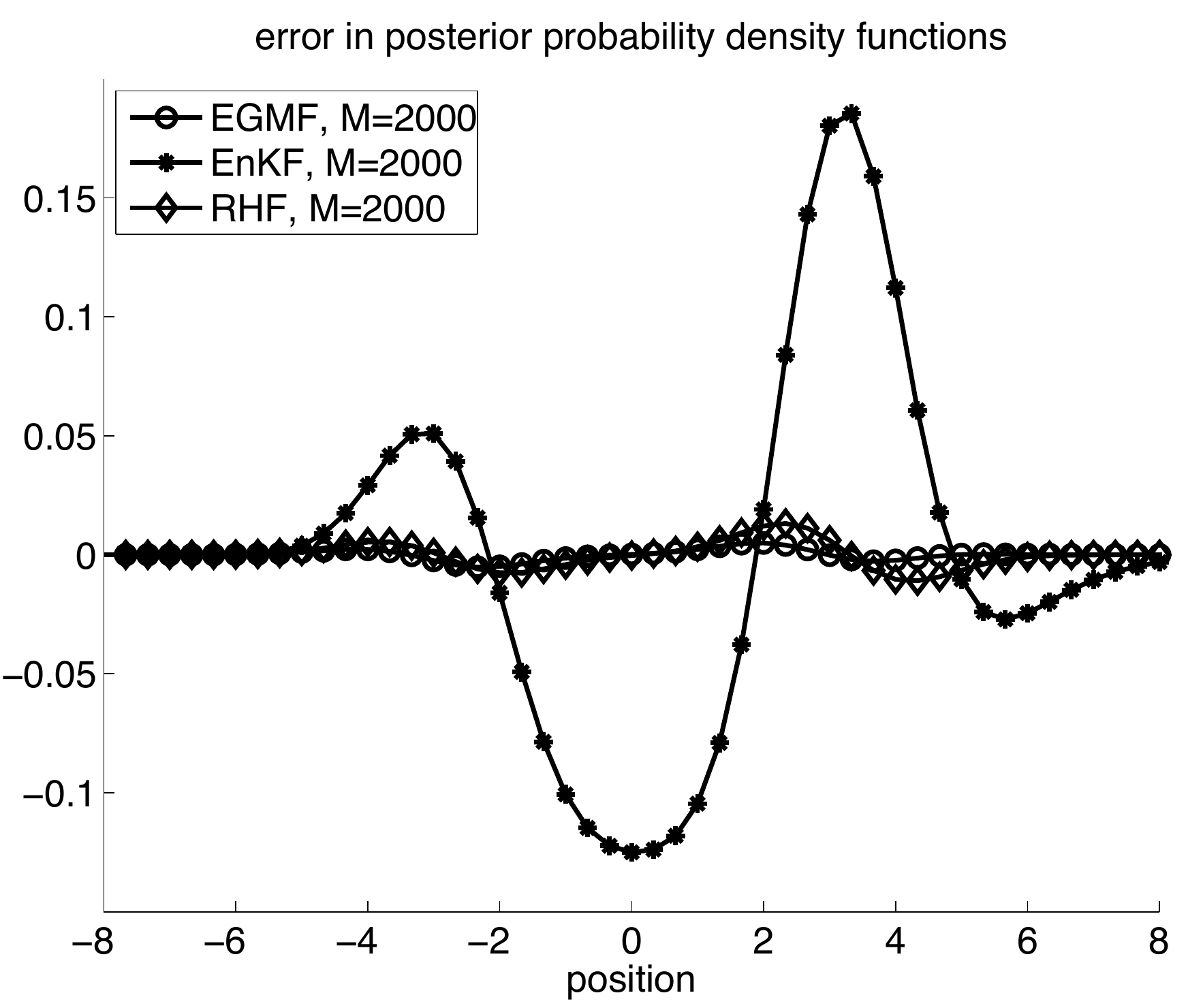}
\end{center}
\caption{Numerically obtained posterior for ensemble sizes $M=50$
  (left panel) and $M=2000$ (right panel). Shown are results from the
EGMF, the RHF, and an EnKF analysis step. While the EGMF and the
RHF converge to the correct posterior distribution, the EnKF 
leads to qualitatively incorrect results for both ensemble sizes.}
\label{figN1.2}
\end{figure}

\begin{figure}
\begin{center}
\includegraphics[width=0.4\textwidth]{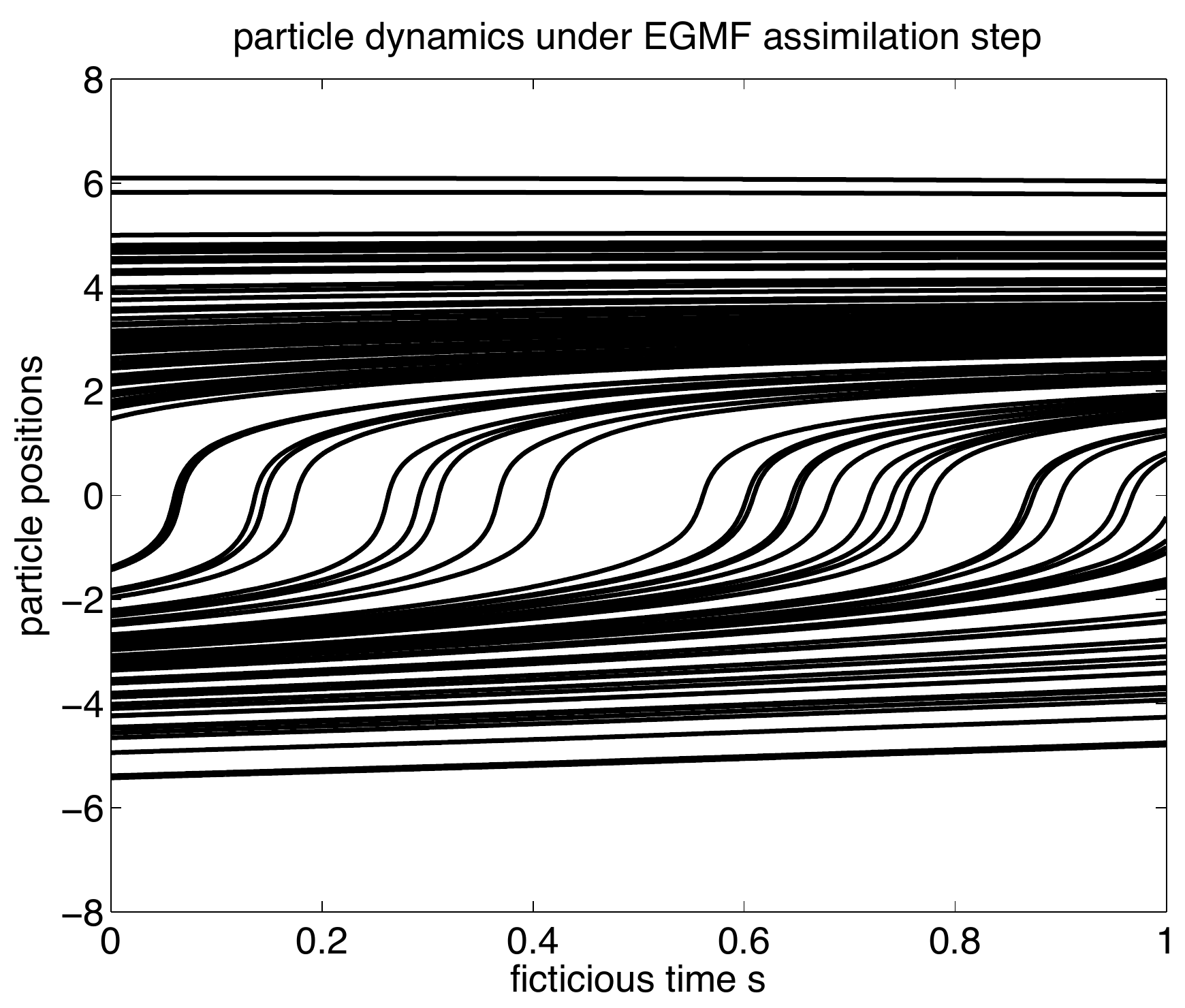}
\end{center}
\caption{Displayed is the rearrangement of the particles under the
  dynamics of the EGMF analysis step. Rapid transitions between the
  Gaussian mixture components are induced by
the vector field $u_{\rm B}$.}
\label{figN1.3}
\end{figure}

We test our formulation first for a single assimilation step where the
prior is a bimodal Gaussian
\begin{equation}
\pi_{\rm prior}(x) = \frac{1}{2} \frac{1}{\sqrt{2\pi}}
e^{-(x-\pi)^2/2} + \frac{1}{2} \frac{1}{\sqrt{2\pi}} e^{-(x+\pi)^2/2}
\end{equation}
and the likelihood is 
\begin{equation}
\pi(y_{\rm obs} | x) = \frac{1}{\sqrt{2\pi}4} e^{-(x-\pi)^2/32} .
\end{equation}
The posterior distribution is again bimodal Gaussian and can be
computed analytically. See Fig.~\ref{figN1.1}. Here we demonstrate how an EnKF, the RHF, 
and the proposed EGMF approximate the posterior for ensemble sizes
$M = 50,2000$ and for $x_i(0) \sim \pi_{\rm prior}$,
$i=1,\ldots,M$. See Fig.~\ref{figN1.2}. Both the RHF and the EGMF are
capable of reproducing the Bayesian assimilation step correctly for $M$
sufficiently large while the EnKF leads to a qualitatively incorrect
result. The transformation of the ensemble members (particles) under the
dynamics (\ref{total}) is displayed in Fig.~\ref{figN1.3}.\\


\begin{figure}
\begin{center}
\includegraphics[width=0.4\textwidth]{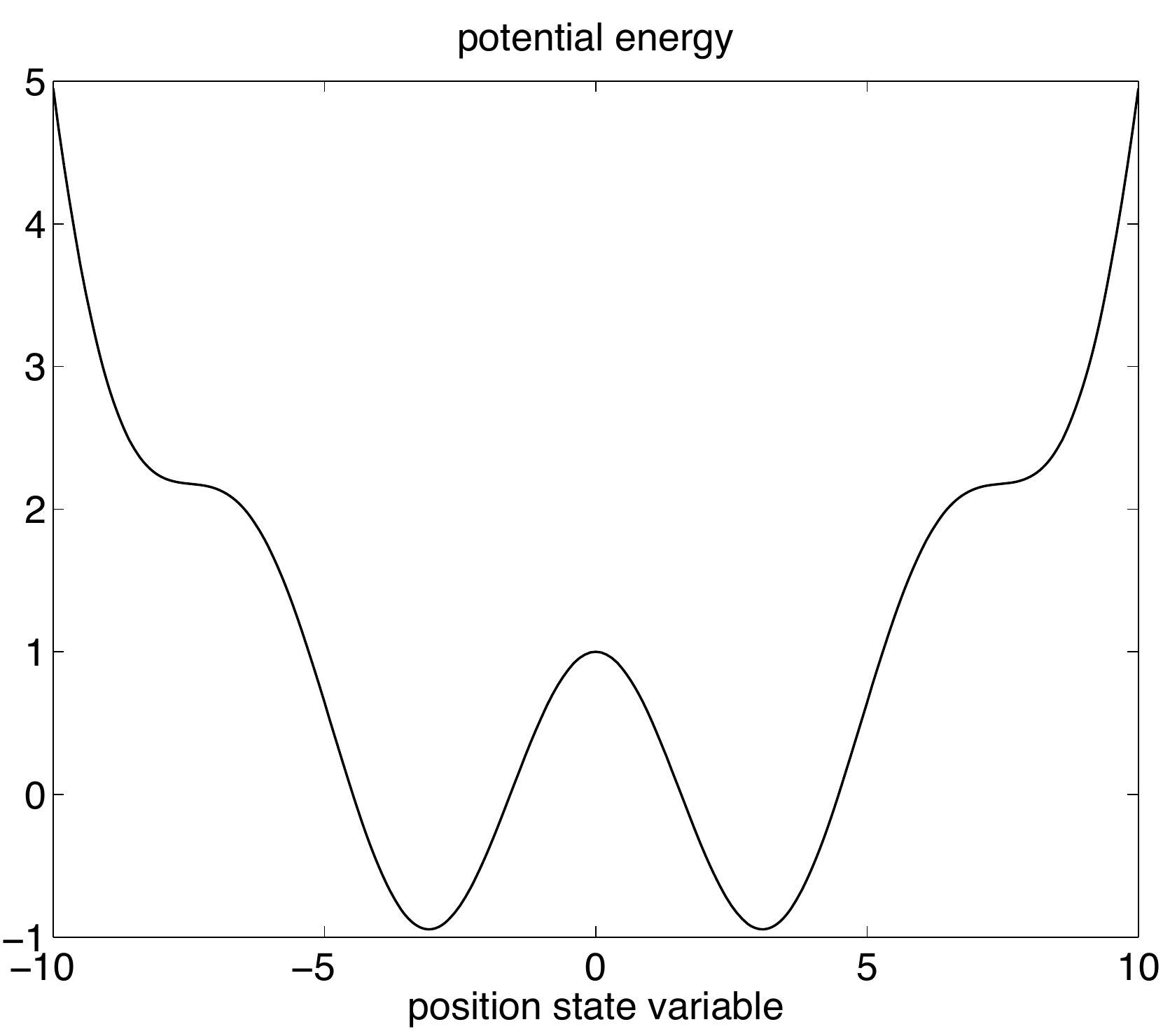} $\qquad$
\end{center}
\caption{Shown is the potential energy $V$ used in both the
  Brownian and Langevin dynamics model.}
\label{figm1}
\end{figure}

\begin{figure}
\begin{center}
\includegraphics[width=0.4\textwidth]{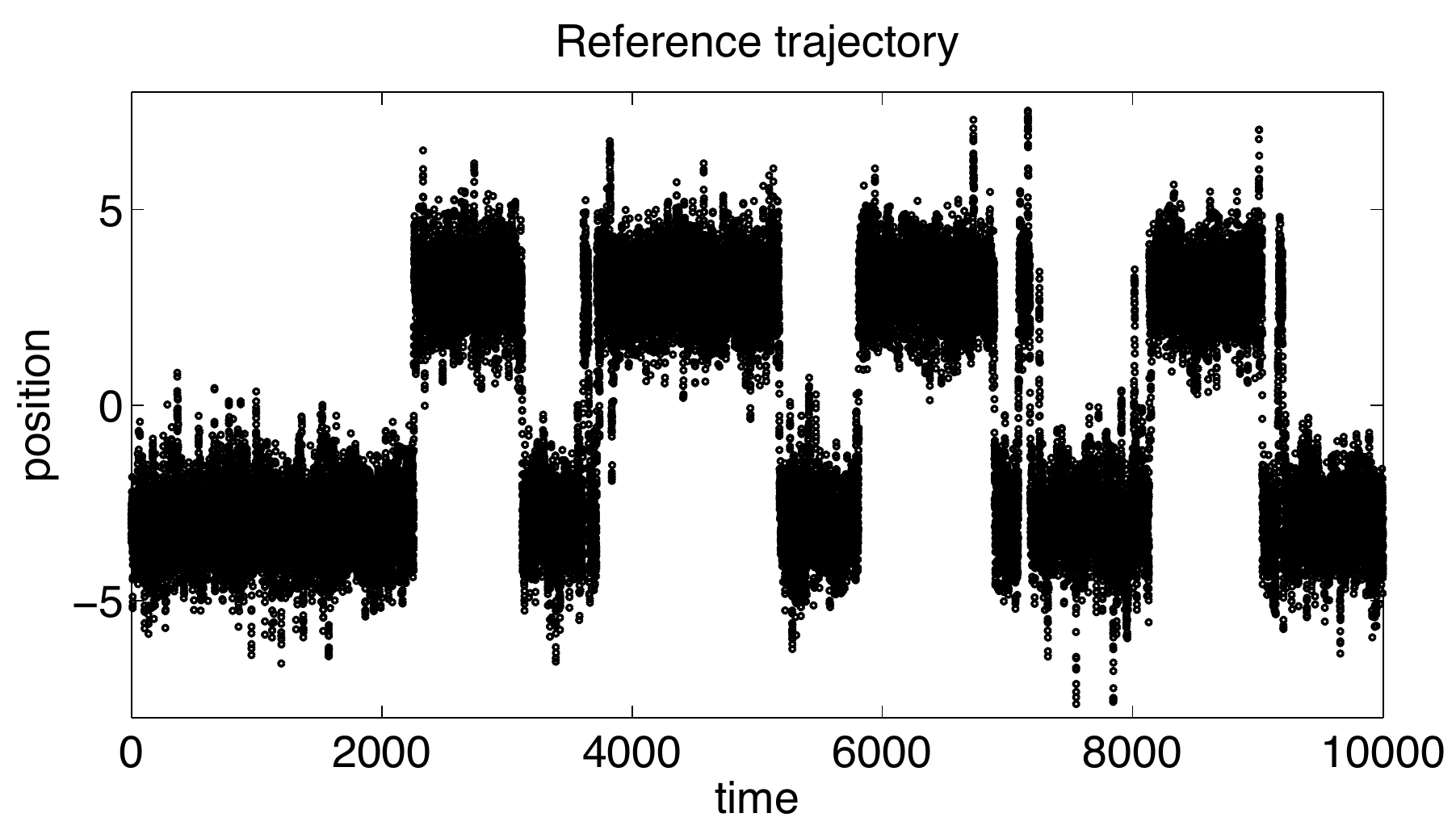}
\end{center}
\caption{Shown is the reference solution from which observations are
  generated by adding Gaussian noise with mean zero and variance $R$.}
\label{fig0}
\end{figure}

\begin{figure}
\begin{center}
\includegraphics[width=0.4\textwidth]{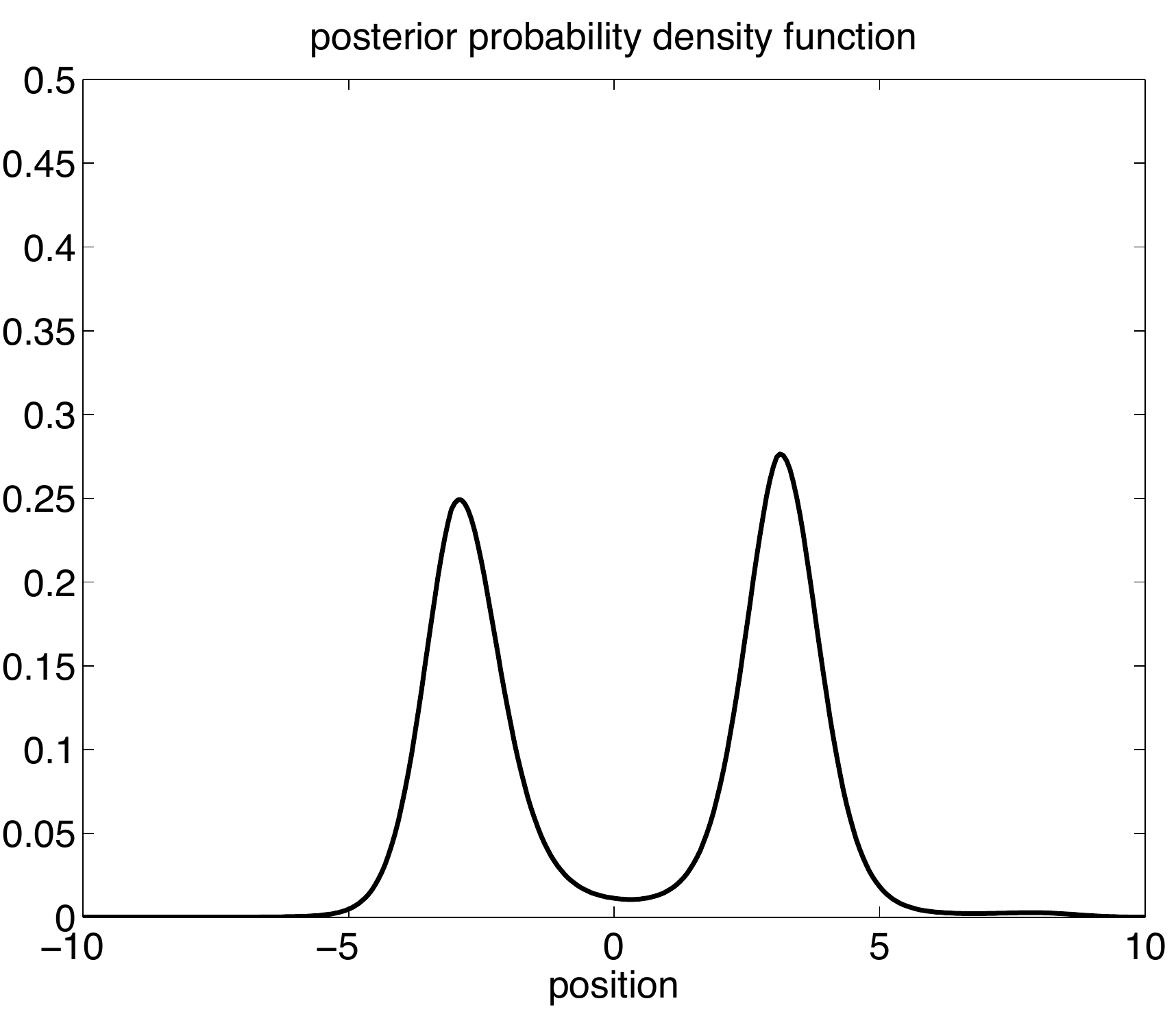}
$\qquad$
\includegraphics[width=0.4\textwidth]{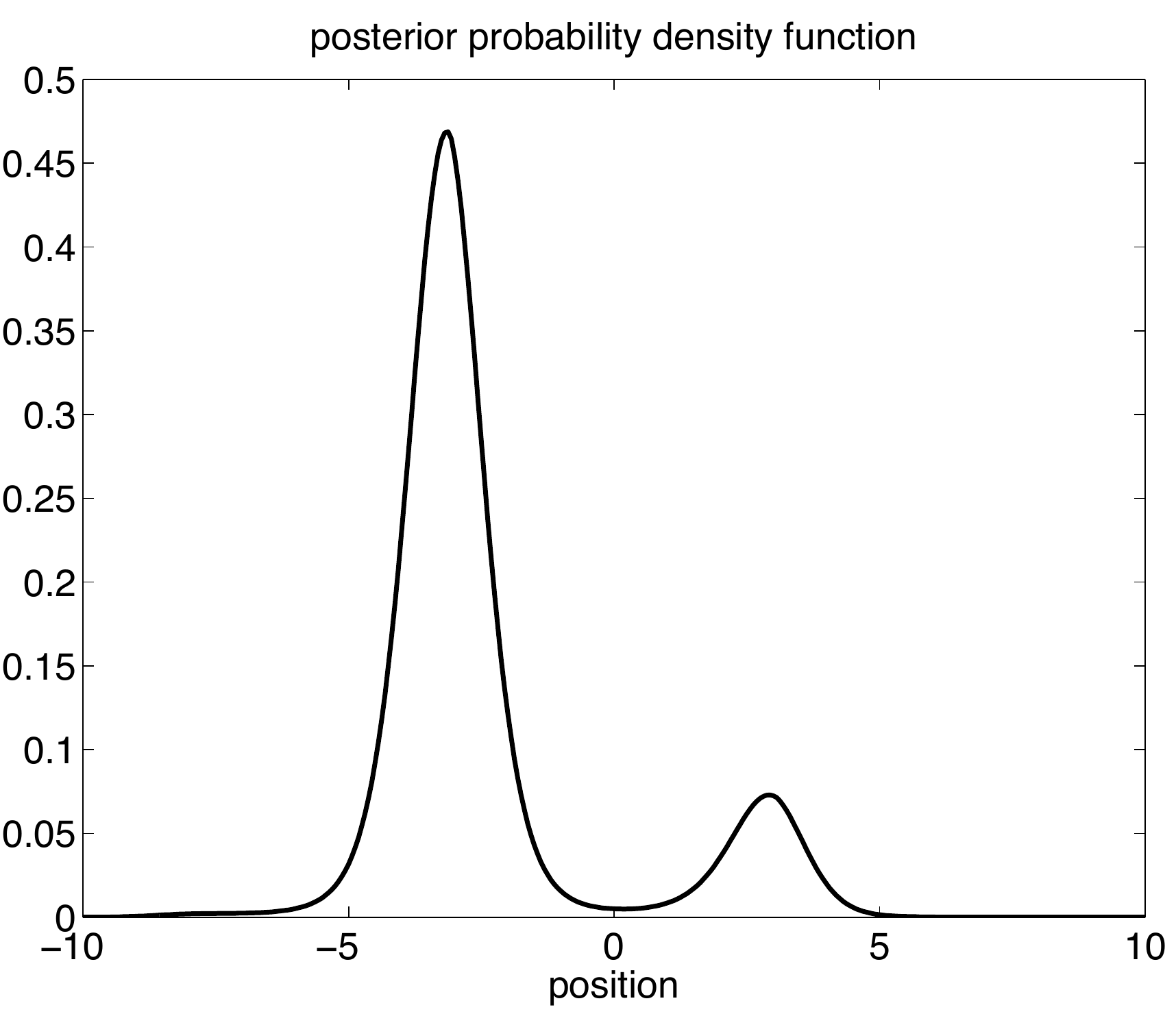}
\end{center}
\caption{Displayed are two posterior PDFs for 
measurement error variance $R=36$ as obtained from the Fokker-Planck
approach. A distinct bimodal behavior can be observed which
motivates the use of a binary Gaussian mixture model for the EGMF.}
\label{fig1}
\end{figure}

\begin{figure}
\begin{center}
\includegraphics[width=0.4\textwidth]{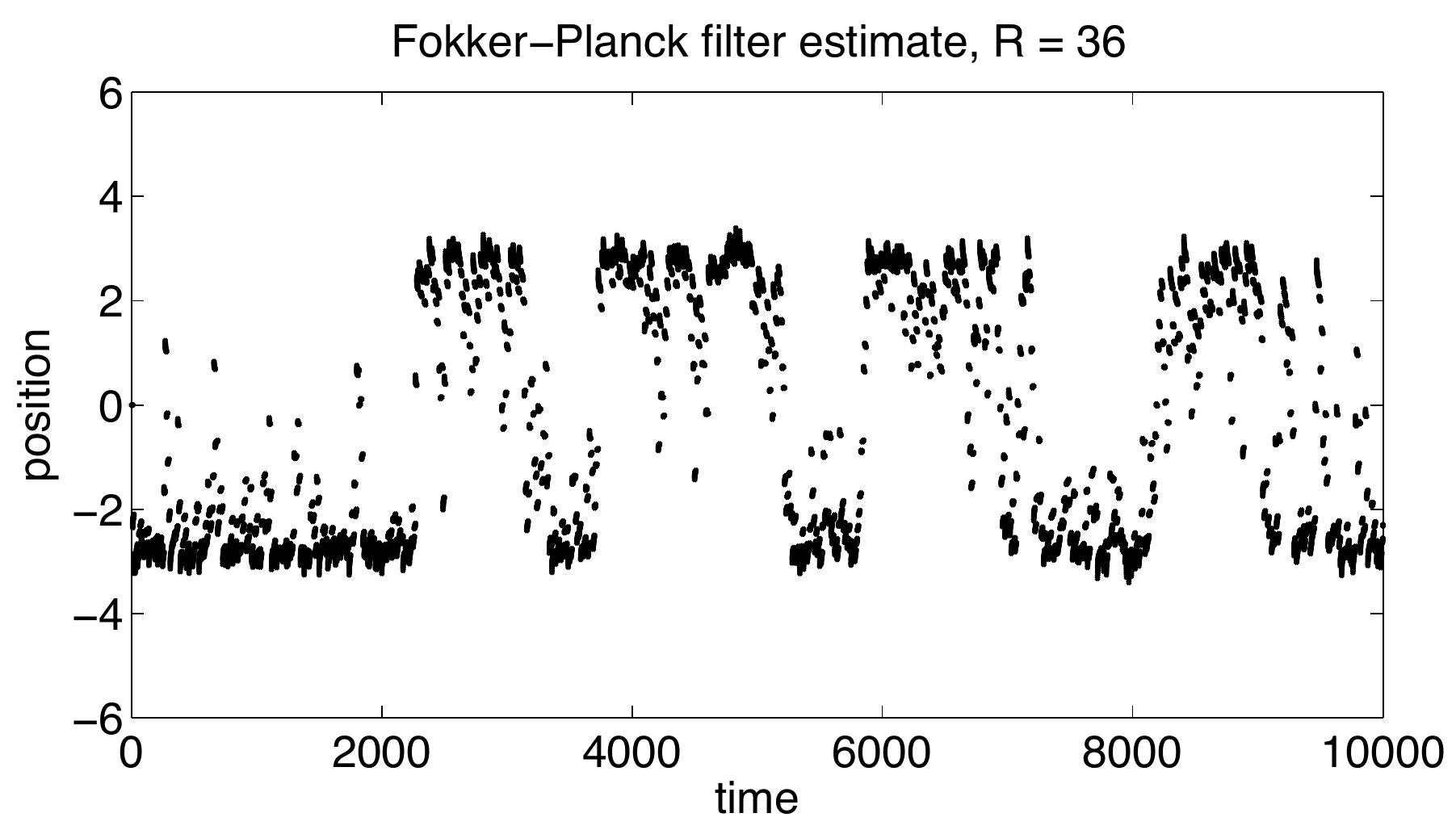} \qquad
\includegraphics[width=0.4\textwidth]{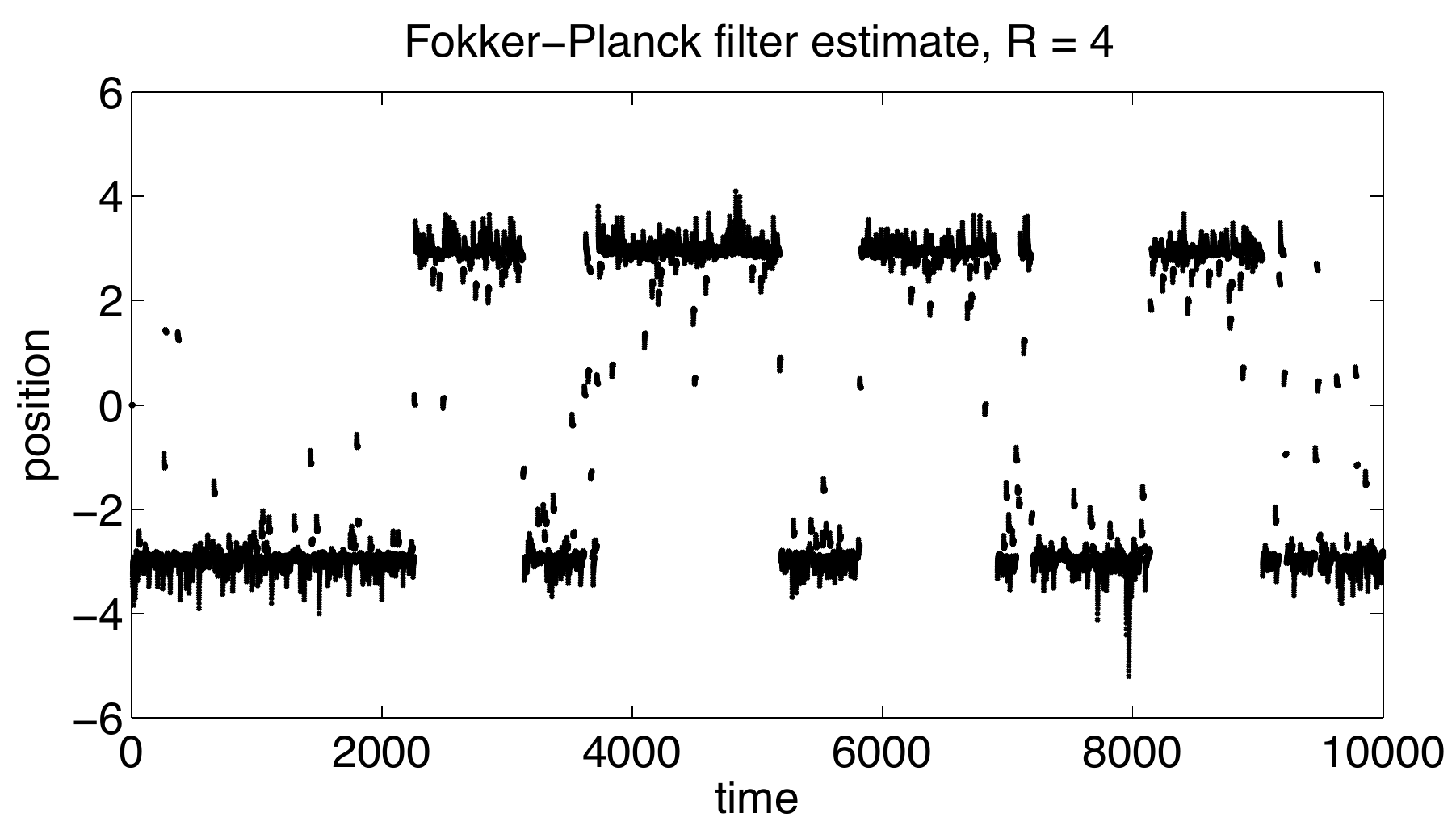}
\end{center}
\caption{Estimated ensemble mean computed from a direct
  simulation of the assimilation problem using a discretized
  Fokker-Planck equation  for
measurement error variance $R=36$ (left panel) and $R=4$ (right panel).}
\label{fig6}
\end{figure}

\begin{figure}
\begin{center}
\includegraphics[width=0.4\textwidth]{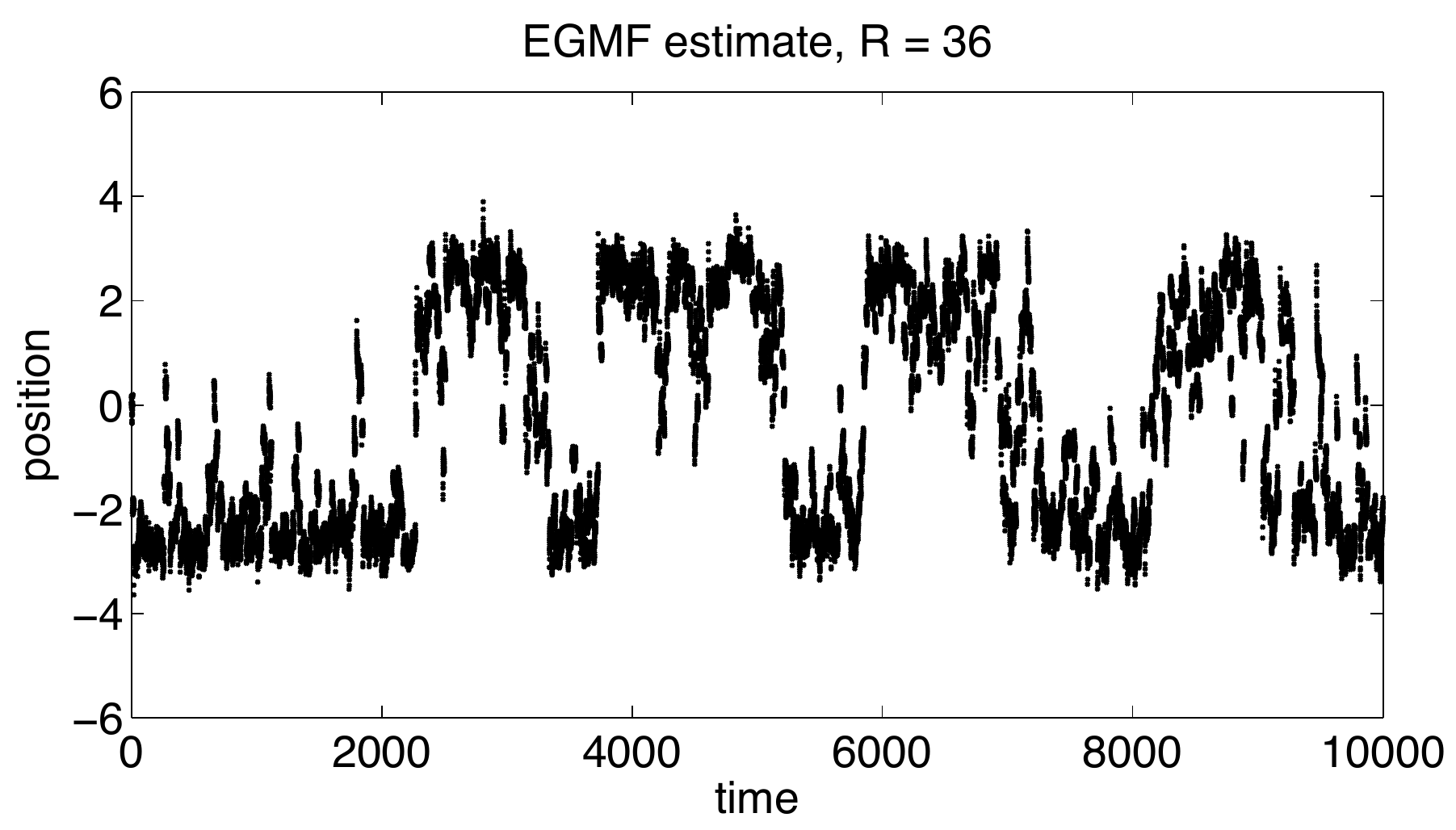}
\qquad
\includegraphics[width=0.4\textwidth]{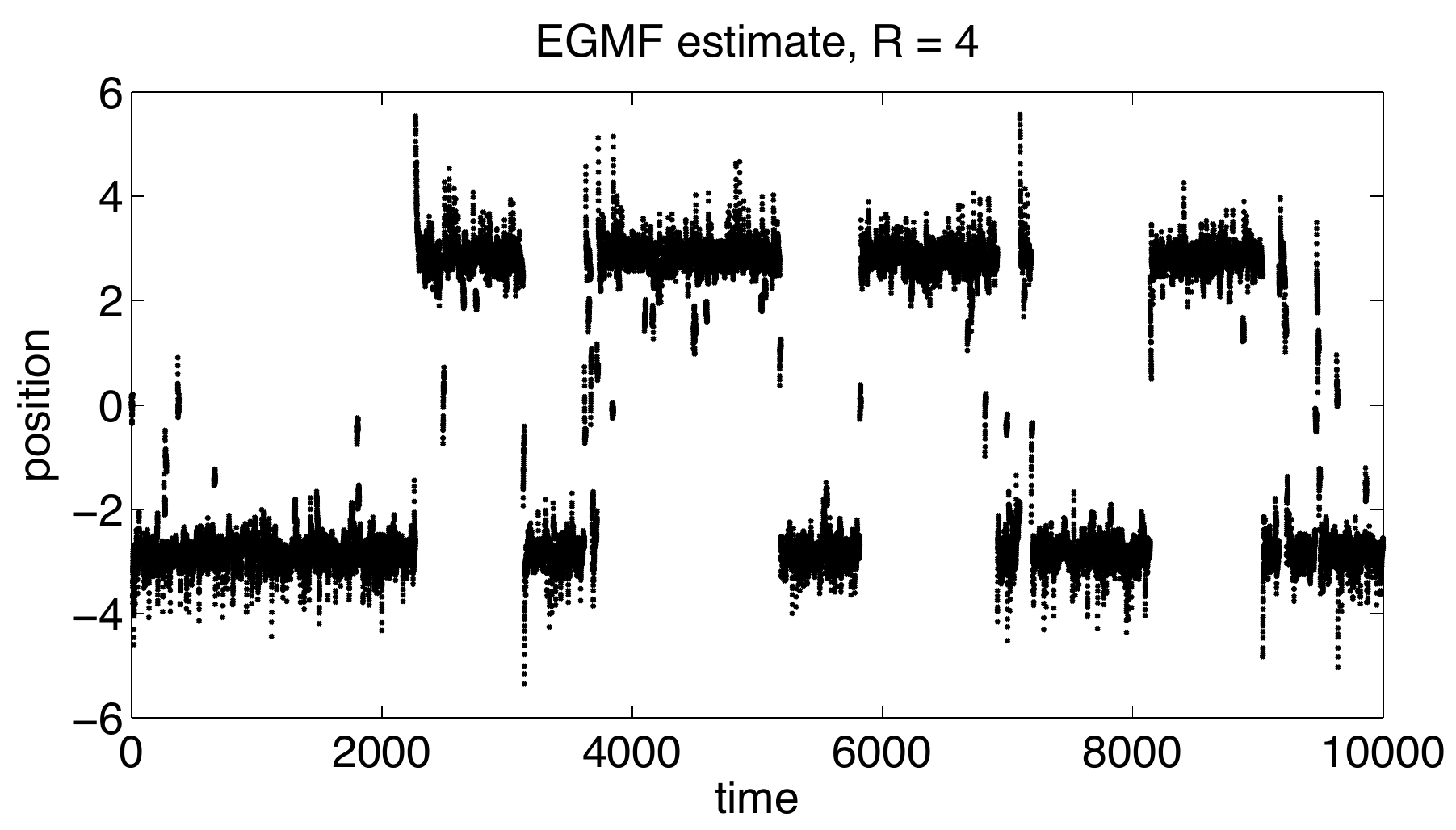}
\end{center}
\caption{Ensemble mean from the EGMF  for
measurement error variance $R=36$ (left panel) and $R=4$ (right panel). It can be
observed that the EGMF leads to results similar to those from the
discrete Fokker-Planck approach  (Fig.~\ref{fig6}).}
\label{fig2}
\end{figure}

\begin{figure}
\begin{center}
\includegraphics[width=0.4\textwidth]{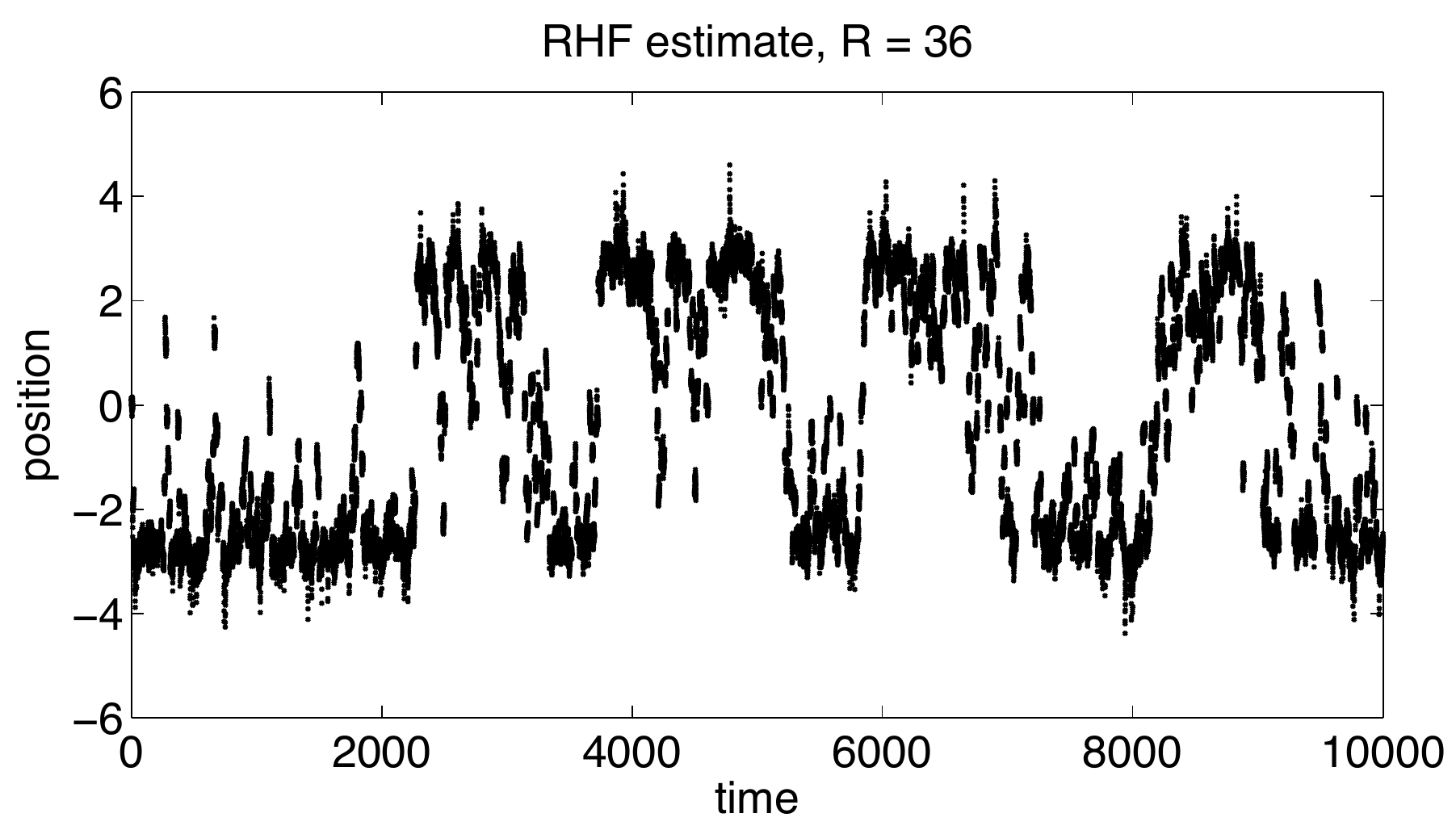}
\qquad
\includegraphics[width=0.4\textwidth]{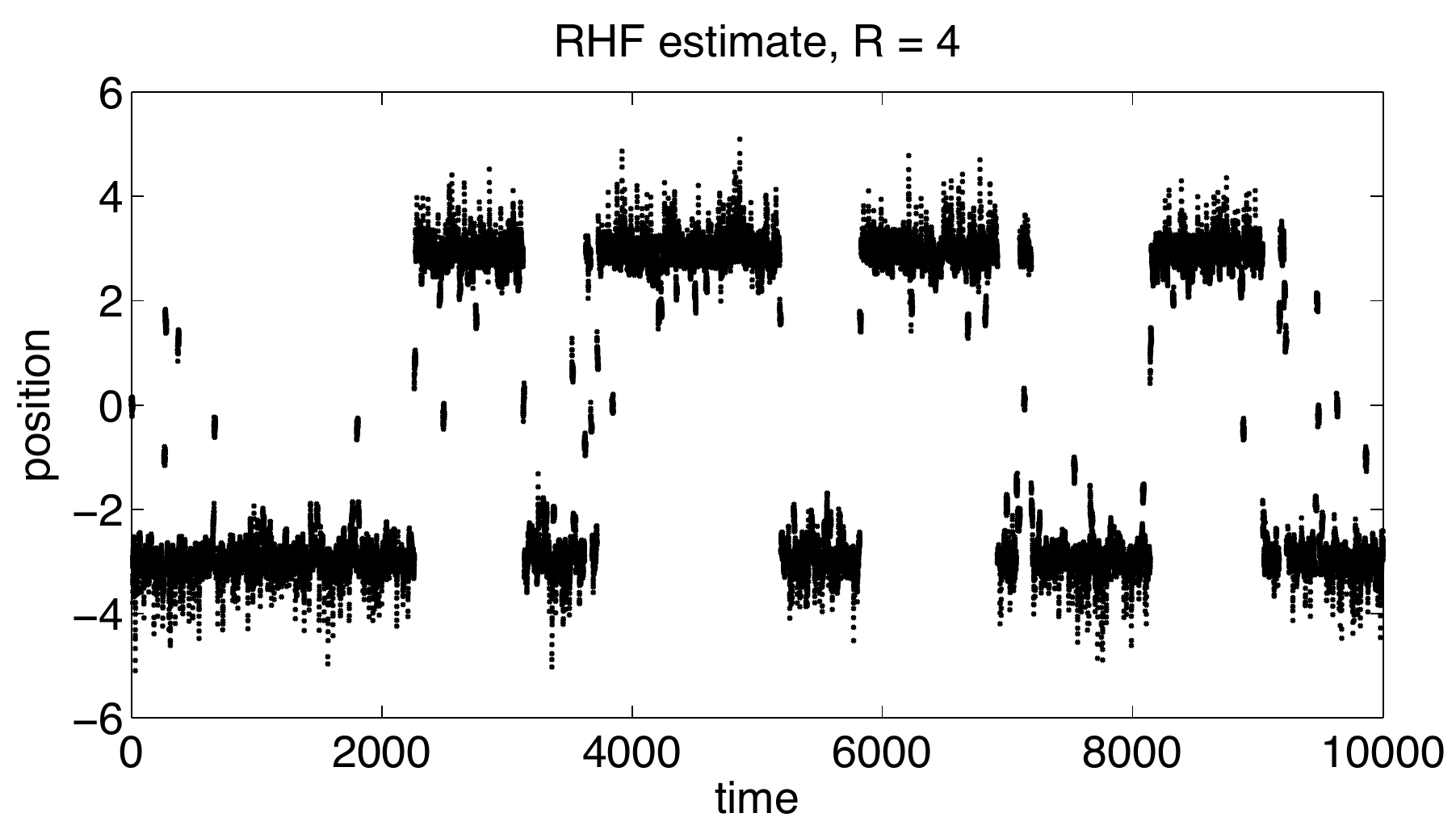}
\end{center}
\caption{Ensemble mean from a RHF
for measurement error variance $R=36$ (left panel) and
  $R=4$ (right panel). The results are for $M=50$ particles. 
It can be observed that the RHF leads to results similar to those from the
discrete Fokker-Planck approach  (Fig.~\ref{fig6}).}
\label{fig3}
\end{figure}

\begin{figure}
\begin{center}
\includegraphics[width=0.4\textwidth]{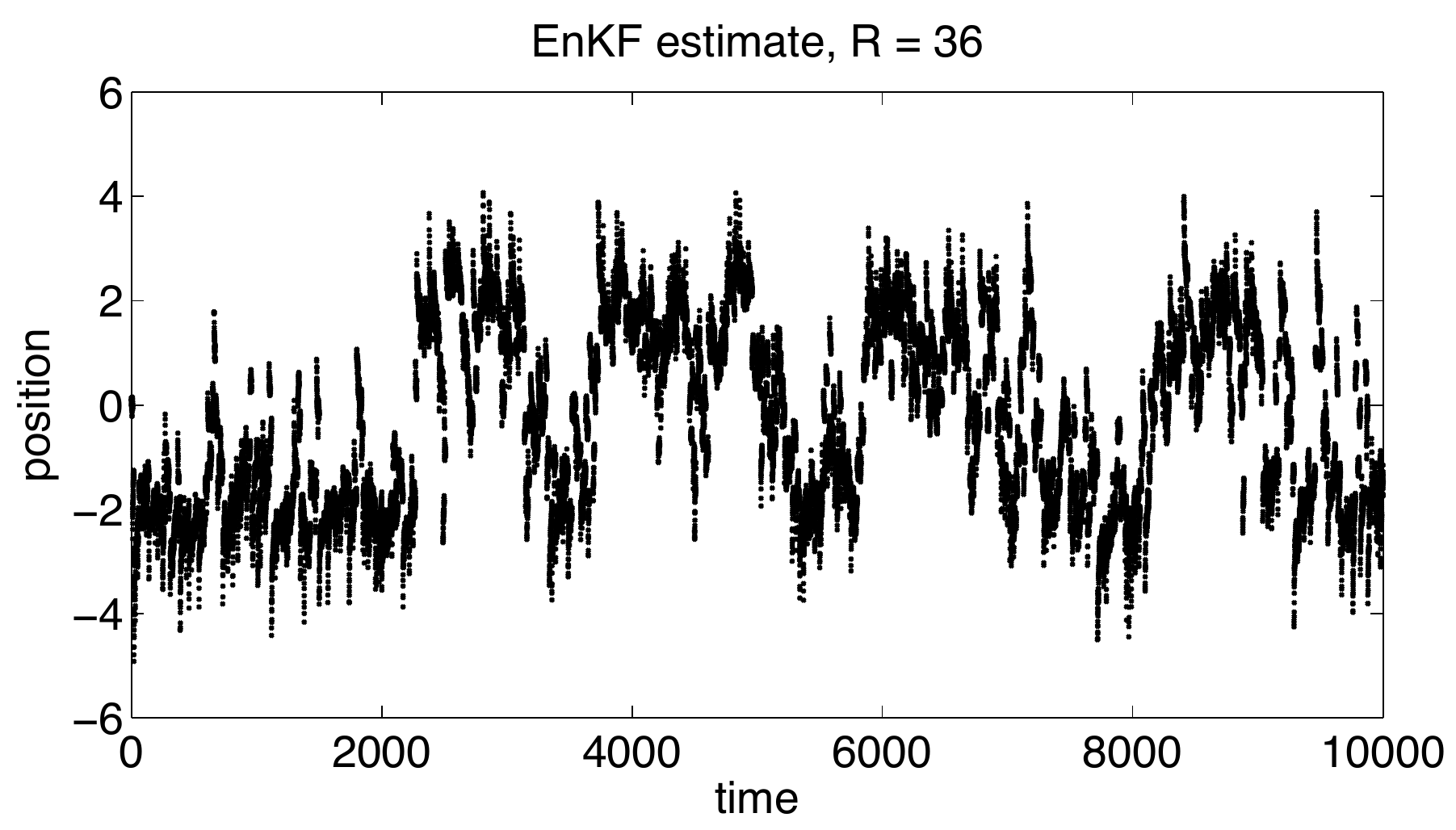}
\qquad
\includegraphics[width=0.4\textwidth]{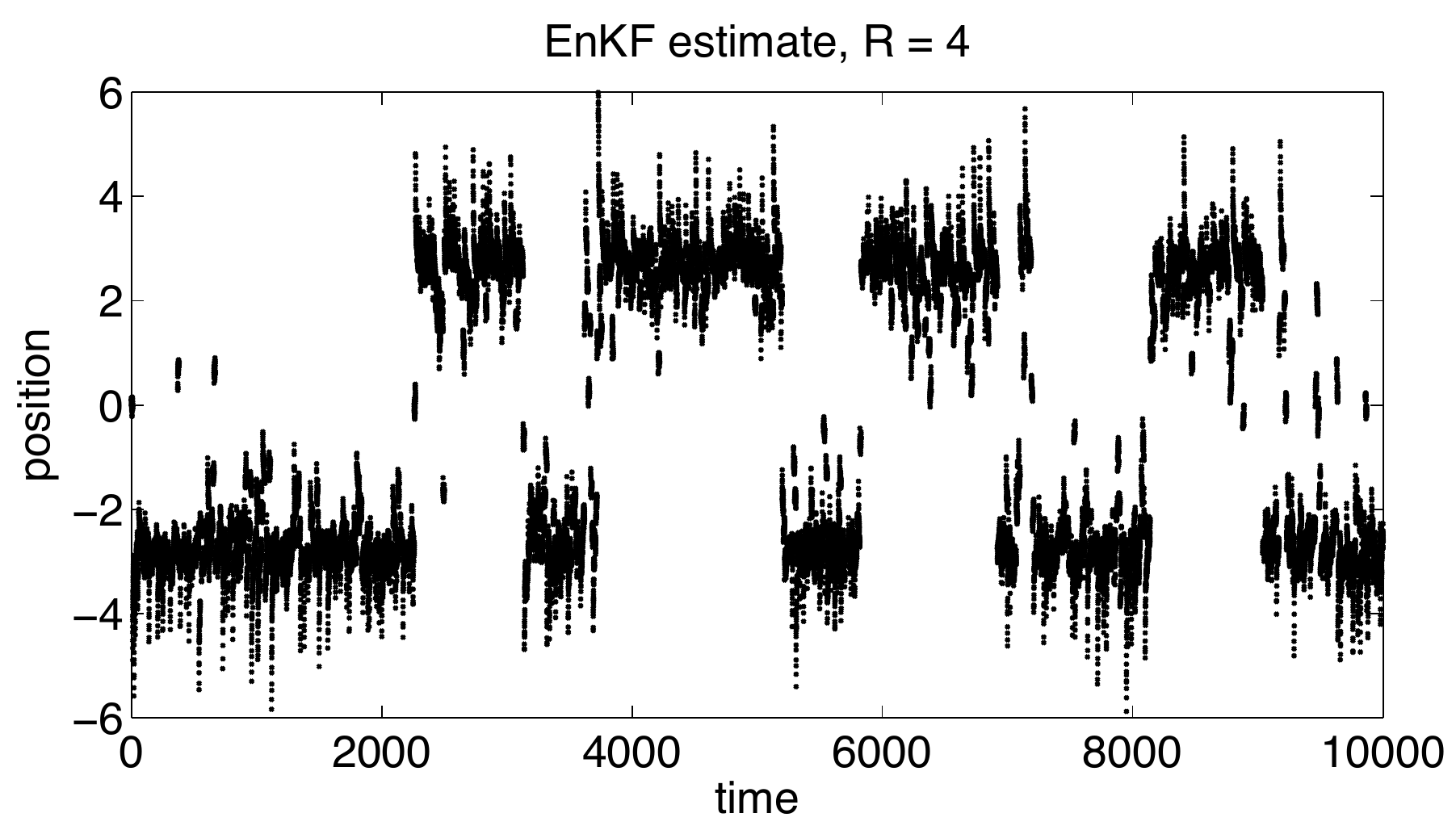}
\end{center}
\caption{Ensemble mean from an EnKF with perturbed
  observations 
 for measurement error variance $R=36$ (left panel) and
  $R=4$ (right panel). The results are for $M=50$ particles. The EnKF
  behaves not as well as the RHF,
  the EGMF, and the discretized Fokker-Planck approach. Similar results are obtained for an
  ensemble square root filter. }
\label{fig4}
\end{figure}

\noindent
We now present results from three increasingly complex filtering problems.


\subsection{A one-dimensional model problem}

As a first numerical example we consider Brownian dynamics under a one-dimensional
potential $V(x)$, i.e.
\begin{equation} \label{NE1}
{\rm d} x = -V'(x){\rm d} t + {\rm d} w(t),
\end{equation}
where $w(t)$ denotes standard Brownian motion and the potential is given by
\begin{equation}
V(x) = \cos(x) + \frac{3}{4}(x/6)^4 .
\end{equation}
See Fig.~\ref{figm1}. The true reference trajectory is started at $x(0) = -3.14$.  See Fig.~\ref{fig0}.
Measurements of $x(t)$ are collected every 10 time units with two different values of the 
measurement error variance $R$ ($R=4,36$).  

The initial PDF is given by the bimodal distribution 
\begin{equation}
\pi_0(x) \propto e^{-(x-3.14)^2/2} + e^{-(x+3.14)^2/2}.
\end{equation}
Depending on the distribution of ensemble members we either use a
single Gaussian ($L=1$) or a bi-Gaussian mixture ($L=2$) in the EGMF
assimilation step. A single Gaussian is used whenever more than 90\%
of the particles are located near either the right (i.e.~$x_i > 0$) or left potential
well (i.e.~$x_i < 0$).  The computed variances are modified such that
$\sigma_l^2  \ge 0.0005$ to avoid singularities in the EM algorithm.
The model equation (\ref{NE1}) is discretized with the forward Euler method
and time-step $\Delta t = 0.1$. The total simulation interval is $t
\in [0,100000]$. The assimilation equations with (\ref{EGM1}) and
(\ref{EGM2}) are discretized with the forward Euler method
and step-size $\Delta s = 0.05$. The $l_\infty$-norm of $u_{\rm B}({\bf
  x},s)$ is limited to a value of $u_{\rm cut} = 5/\Delta s$. 

The performance of the EGMF
is compared to an EnKF with perturbed observations, 
an ensemble square root filter, and a RHF. 
The particle positions are adjusted during each data assimilation step of a RHF 
such that the particles maintain equal weights $\gamma_i = 1/M$. The adjustment is
done similar to what has been proposed by \cite{sr:anderson10} except
that the posterior is approximated by piecewise constant functions.

For this simple, one-dimensional problem the densities can be directly
propagated through a discretization of the associated Fokker-Planck
equation. Bayes theorem reduces to a multiplication of the prior
PDF approximation from the Fokker-Planck approximation with the
likelihood at each grid point. We have used a grid with mesh-size
$\Delta x = 0.125$ over $x\in [-10,10]$ to provide an independent and
accurate filtering result. Periodic boundary conditions are used for the diffusion operator such that
the spatial discretization leads to a stochastic matrix. It is found from the numerical simulations 
that $R=36$ leads to a pronounced bimodal behavior of the
solution PDF $\pi$. See Figure \ref{fig1} for two posterior
PDF approximations from the Fokker-Planck approach.

Typical filter behaviors over the time interval $t\in [0,10000]$
with regard to the reference trajectory can be found in
Figures \ref{fig6}, \ref{fig2}, \ref{fig3}, and \ref{fig4}, respectively, for $M=50$ ensemble 
members. The EGMF and the RHF display a behavior similar to that from the
discretized Fokker-Planck approach while
significantly different results are obtained from the EnKF 
implementation with perturbed observations. Similar results are obtained for 
an ensemble square root filter (not displayed). The EGMF uses a bi-Gaussian approximation in 97\%
of the assimilation steps for $R=36$ and in 47\% of the assimilation
steps for $R=4$.

We also provide the root mean square (RMS) error between the 
computed mean from the three different filters and the mean computed from the
Fokker-Planck approach in Table \ref{table1} for $R=36$ and
different ensemble sizes $M$.  The RHF
converges as $M\to \infty$ to the solution of the Fokker-Planck
approach for this one-dimensional model problem. 
The EGMF provides better results than the EnKF but does not
converge since the limiting distributions are not exactly
bi-Gaussian.  Note that the EGMF should converge for $M\to \infty$ and
the number of mixture components sufficiently large. Note also that
the RHF does not converge to the analytic filter solution as $M\to
\infty$ in case of more than one dynamic variable (i.e.~$N>1$). See also
the following example.

\begin{table}
\caption{RMS errors for ensemble means obtained from EnKF, RHF, and EGMF 
compared to the expected value computed by  a Fokker-Planck
discretization with error variance $R=36$ and $M=20,50,100$ 
particles/ensemble members.}
\label{table1}
\begin{center}
\begin{tabular}{cccc}
\hline\noalign{\smallskip}
& RHF & EGMF  & EnKF\\
\noalign{\smallskip}\hline\noalign{\smallskip} 
$M= 20$ & 0.6551 & 0.7683 & 1.0283 \\
$M= 50$& 0.3717 & 0.5127 & 0.8798 \\ 
$M = 100$ & 0.2691 & 0.4033 & 0.8412 \\ 
\noalign{\smallskip}\hline
\end{tabular}
\end{center} 
\end{table}


\subsection{A two-dimensional model problem}

We discuss another example from classical mechanics. The evolution of
a particle with position $q\in \mathbb{R}$ and velocity $v\in \mathbb{R}$
is described by Langevin dynamics \citep{sr:gardiner} with equations of motion
\begin{eqnarray}
dq &=& v\, dt,\\
dv &=& -V'(q)\,dt - \gamma v\,dt + \sigma dw(t) ,
\end{eqnarray}
where the potential $V(q)$ is given by
\begin{equation}
V(q)= \cos(q) + \frac{3}{4} (q/6)^4 ,
\end{equation}
the friction coefficient is $\gamma = 0.25$, $w(t)$ denotes standard
Brownian motion, and $\sigma^2 = 0.35$. A reference solution, denoted
by $(q_r(t),v_r(t))$, is obtained for initial condition $(q_0,v_0) = (1,1)$ and a particular
realization of $w(t)$. 

Let us address the situation that 
the reference solution is not directly accessible to us and that instead we are only
able to observe $Q(t)$ subject to
\begin{equation} \label{obs}
dQ(t) = v_r(t)\,dt + c^{1/2} d\xi(t),
\end{equation}
where $\xi(t)$ denotes again standard Brownian motion and $c =
0.2$. In other words, we are effectively only able to observe particle
velocities. 

We now combine the model equations and the
observations within the ensemble Kalman-Bucy framework. The ensemble
filter relies on the simultaneous propagation of an ensemble of
solutions $x_i(t) = (q_i(t),v_i(t)) \in \mathbb{R}^2$, $i=1,\ldots,M$.  In our experiment we set
$M=50$. The filter equations for an EGMF with a single Gaussian, ensemble Kalman-Bucy filter
\citep{sr:br11}, respectively, become
\begin{eqnarray}
dq_i &=& v_i\,dt - \frac{P_{qv}}{2c} (v_i\,dt +  \bar v \,dt - 2dQ(t)  ),\\
dv_i &=& -V'(q_i)\,dt - \gamma v_i\,dt + \sigma dw_i(t) -
\frac{P_{vv}}{2c} (v_i\,dt +\bar v \,dt
- 2dQ(t) )
\end{eqnarray}
with mean
\begin{equation}
\bar v = \frac{1}{M} \sum_{i=1}^M v_i, \qquad \bar q = \frac{1}{M} \sum_{i=1}^M
q_i
\end{equation}
and variance/covariance
\begin{equation}
P_{vv} = \frac{1}{M-1} \sum_{i=1}^M (v_i-\bar v)^2, \qquad
P_{qv} = \frac{1}{M-1} \sum_{i=}^M  (q_i-\bar q)(v_i-\bar v) .
\end{equation}
The equations are solved for each ensemble member with different
realizations $w_i(t)$  of standard Brownian motion and step-size
$\Delta t = 0.01$. The observation interval in (\ref{obs}) is $\tau =
\Delta t$. The extension of the EGMF to a Gaussian mixture with $L=2$ is
straightforward. One substitutes $y=v$, $R = c/\Delta t$, and $y_{\rm obs} = \Delta Q(t_n)/\Delta t$
with
\begin{equation}
\Delta Q(t_n) = v_r(t_n)\Delta t + \sqrt{c\Delta t} \xi_n, \qquad \xi_n \sim N(0,1),
\end{equation}
into (\ref{EGM1}) and (\ref{EGM2}) and sets $\Delta s = 1$ in the numerical time-stepping procedure
for the assimilation step. The $l_\infty$-norm of $u_{\rm B}({\bf
  x},s)$ is limited to a value of $u_{\rm cut} = 0.25/\Delta s$. 
Assimilation is performed at every model time-step. 
We perform a total of two million time-steps/data assimilation cycles. In the same manner one
can implement a RHF for this problem. 

\begin{figure}[htb]
\begin{center}
\includegraphics[width=0.4\textwidth]{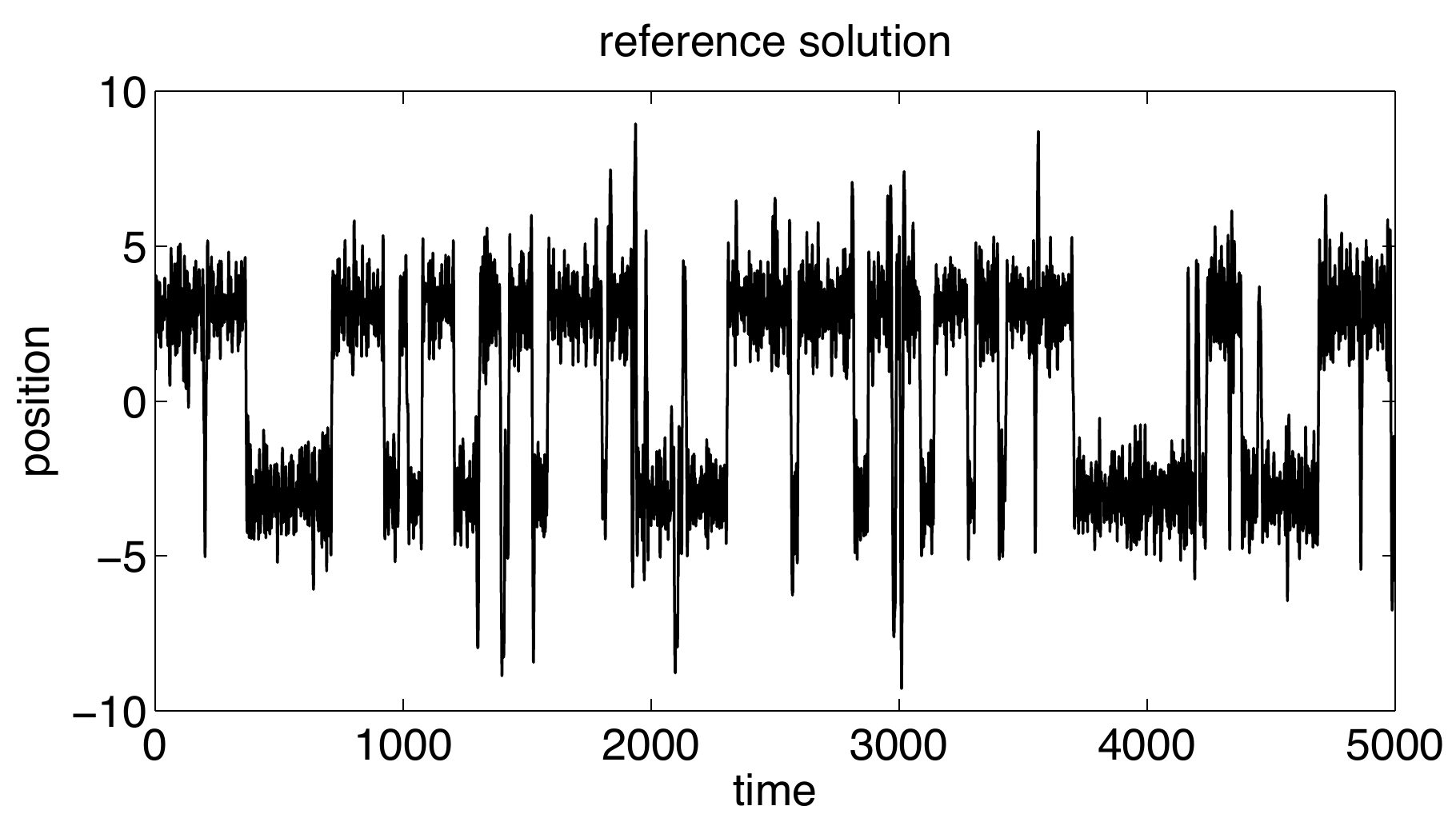} \qquad 
\includegraphics[width=0.4\textwidth]{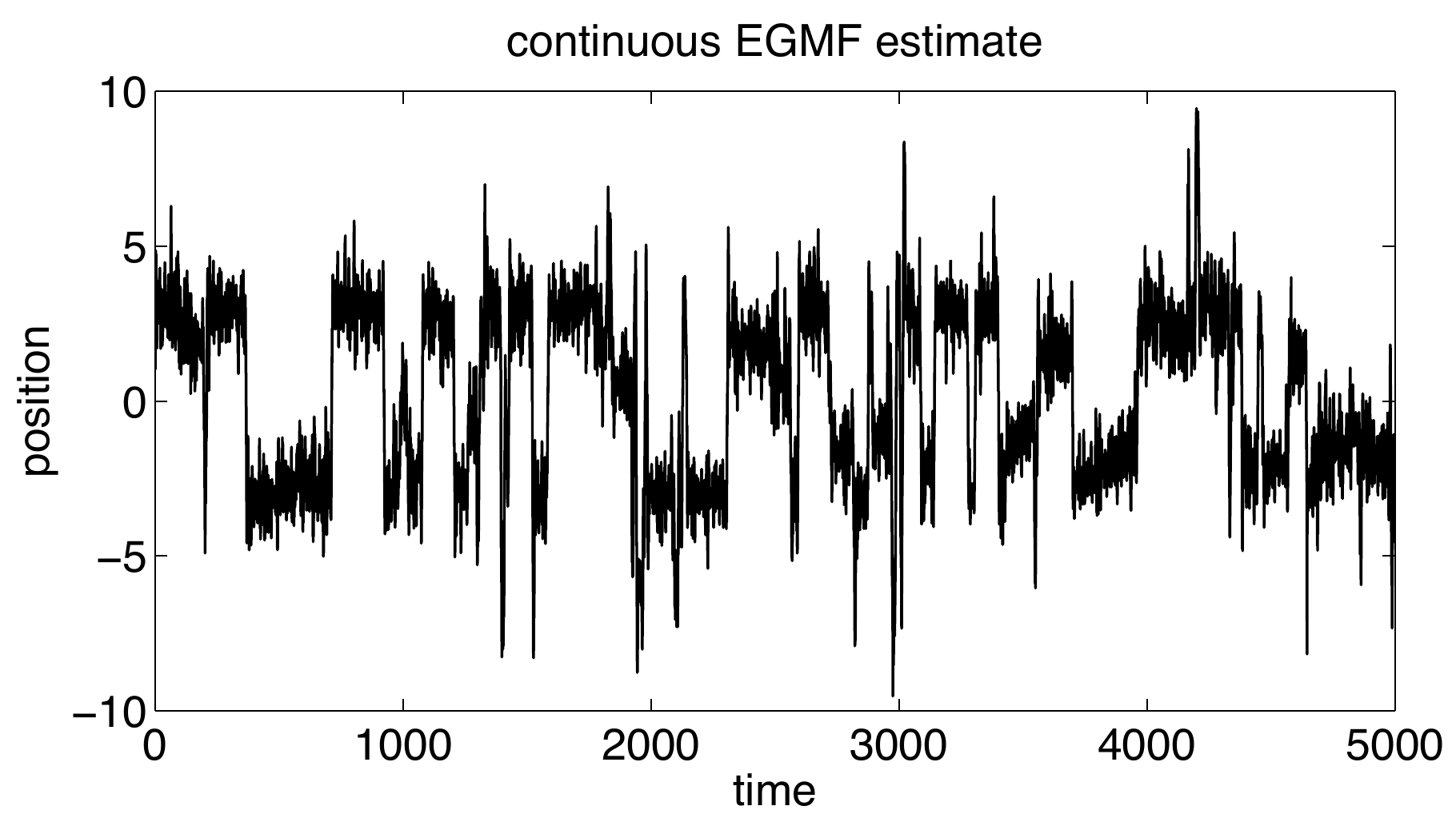} 
\end{center}
\caption{The reference solution $q_r(t)$ (left panel) and the estimated
 (ensemble mean) solution
 from the continuous 
EGMF (right panel) over the first quarter of the simulation interval
$t\in [0,20000]$. The estimated solution mostly follows the reference
solution with the exception of a number of missed transitions.}
\label{figN3.1}
\end{figure}

\begin{figure}[htb]
\begin{center}
\includegraphics[width=0.4\textwidth]{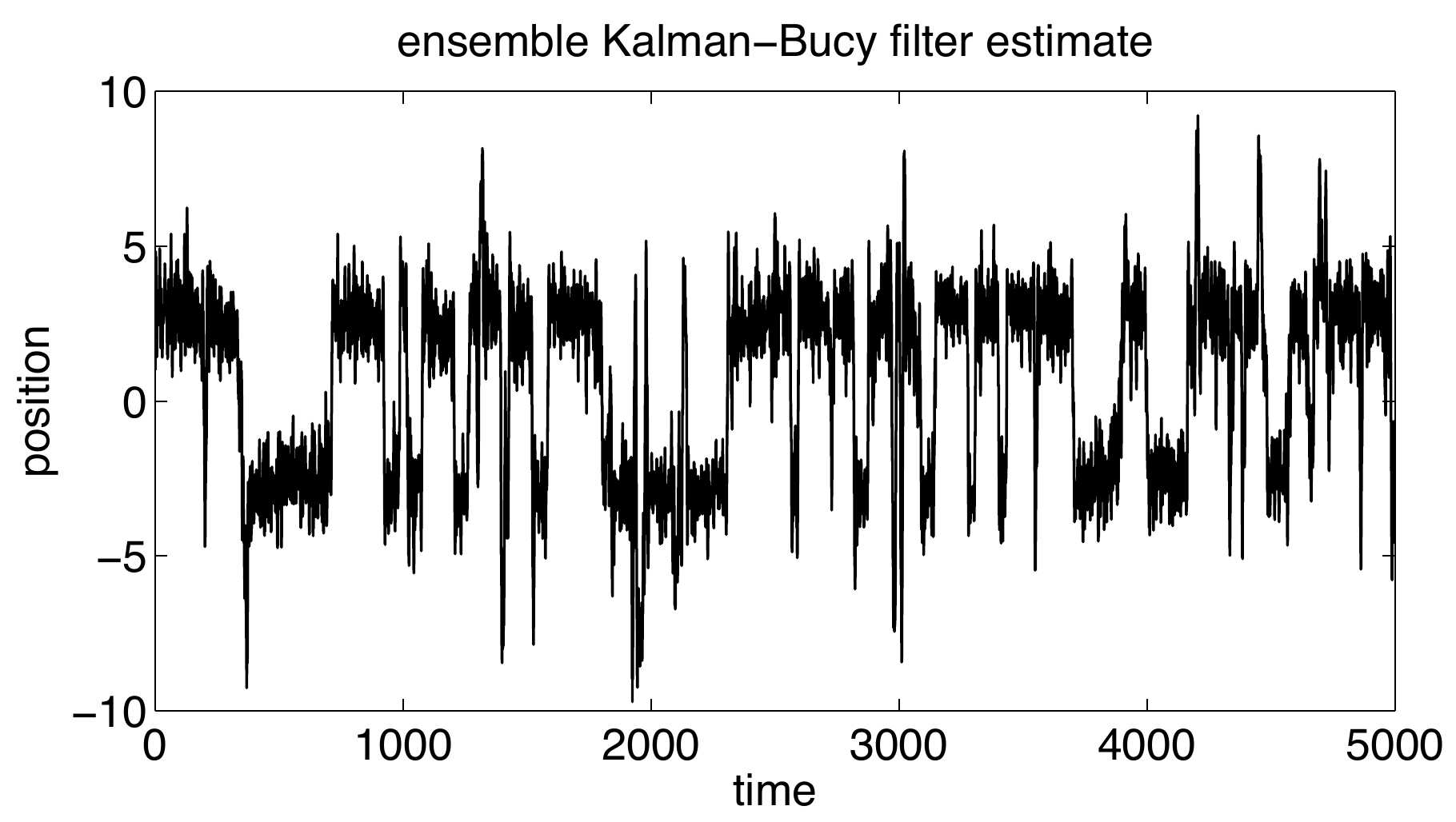} \qquad
\includegraphics[width=0.4\textwidth]{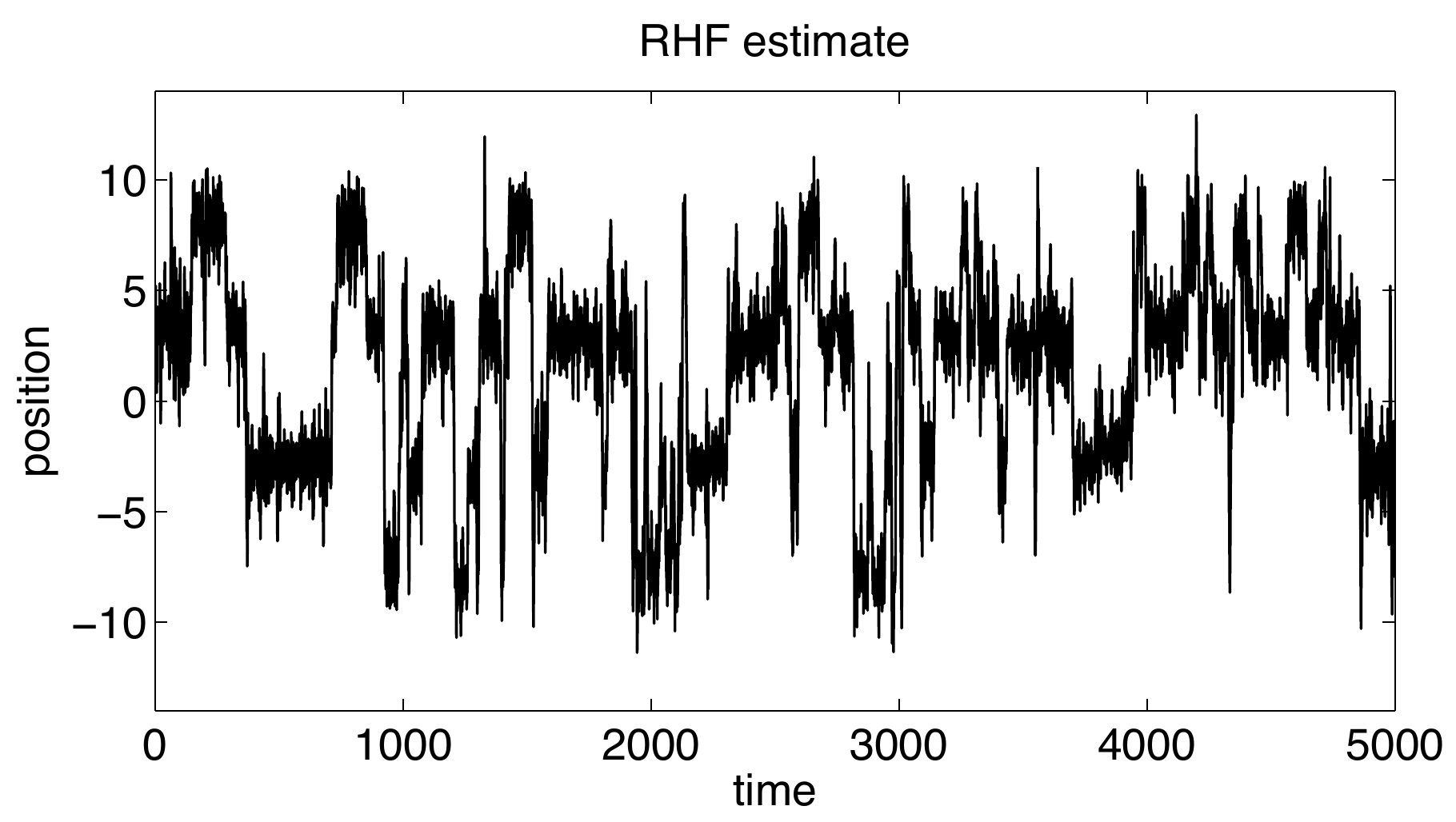} 
\end{center}
\caption{The estimated (ensemble mean) solution from
the ensemble Kalman-Bucy filter over a quarter of the
simulation interval is displayed in the left panel.
The results look qualitatively similar to the results from the EGMF filter. 
However, in terms of RMS errors, the EGMF outperforms the ensemble Kalman-Bucy filter
by about 20\% over the complete simulation interval. 
The right panel displays the estimated (ensemble mean) solution
 from the RHF (left panel). The reader should note the enlarged range of the
 vertical axis.  At several instances the filtered solution strongly deviates from the reference solution. 
 To interpret this behavior we need to have a closer look at the potential energy $V(q)$
 (compare Fig.~\ref{figm1}).  Apparently the RHF interprets the data as
 corresponding to solutions with positions in the flat side 
 regions of the potential energy function.}
\label{figN3.1b}
\end{figure}

The computed ensemble means $\bar
q(t)$ and $\bar v(t)$ in comparison to the reference solution can be
found in Fig.~\ref{figN3.1} for the continuous EGMF (using $L=1$ and $L=2$ mixture 
components as appropriate) and the ensemble Kalman-Bucy filter 
(continuous EGMF with $L=1$). 
The root mean square (RMS) error with respect to the true reference solution 
is 2.3331 for the ensemble Kalman-Bucy filter and 1.9148 for the EGMF, which
amounts to a relative improvement of about 20\%. The EGMF uses a
bi-Gaussian distribution in about 25\% of the assimilation steps. For comparison
we show the results from the RHF in Fig.~\ref{figN3.1b} for $M=50$ particles. To interpret the behavior
of the RHF one needs to look at the potential energy function $V(q)$
(see Fig.~\ref{figm1}). The RHF
assimilation scheme apparently pushes the solutions occasionally into the flat side regions of the
potential energy curve resulting in a relatively large RMS error of 3.9375. Qualitatively similar results
are obtained for the RHF with $M = 800$ particles. Recall that 
we do observe velocities and not positions in this example and that the RHF
uses the ensemble generated covariance matrix ${\bf P}$ to linearly
regress filter increments onto state space.

\subsection{Lorenz-63 model} \label{sec_num3}

\begin{figure}[htb]
\begin{center}
\includegraphics[width=0.4\textwidth]{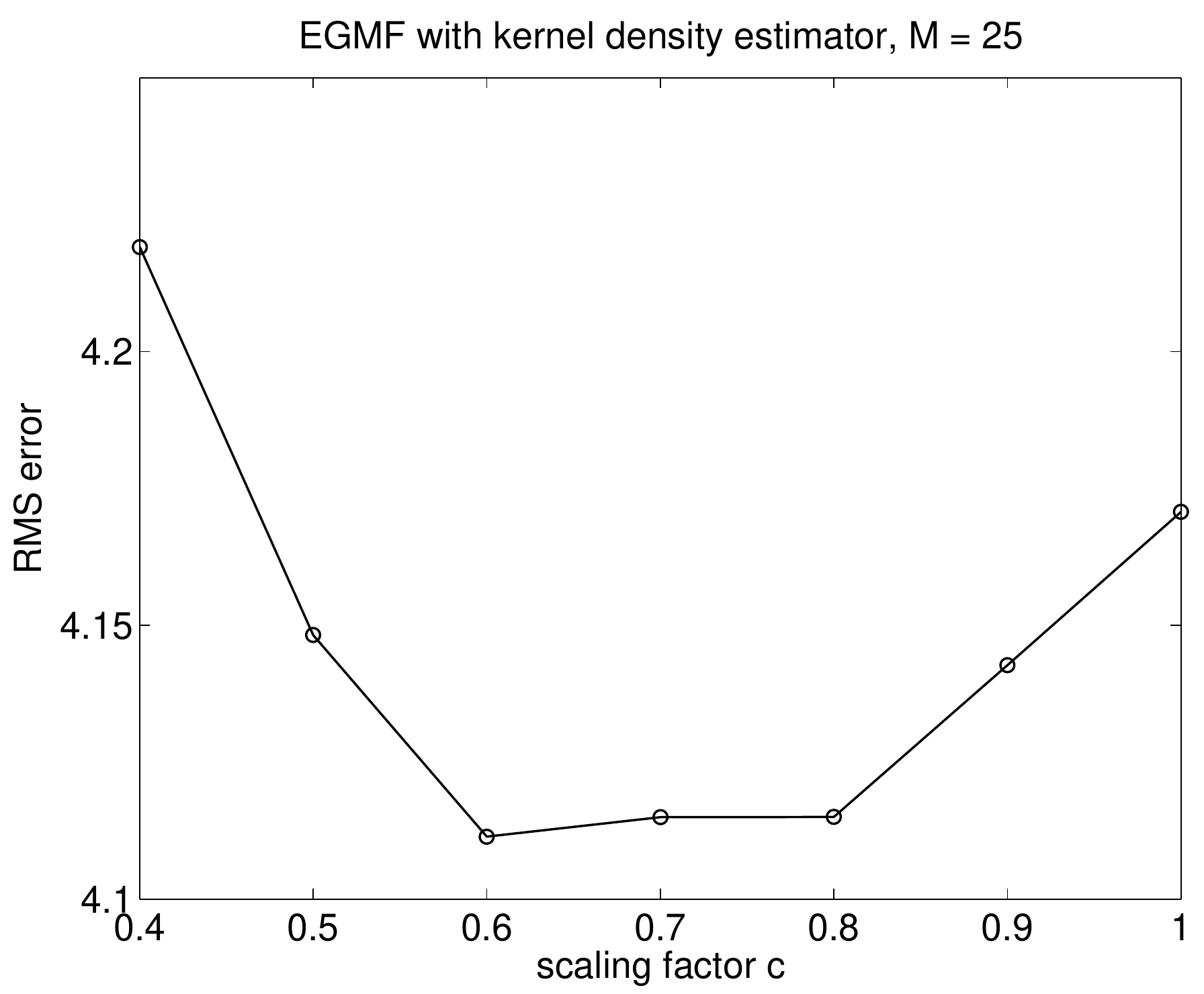}
\end{center}
\caption{RMS errors from the EGMF for a range of values of the scaling
parameter $c$ and ensemble size $M=25$. A series of experiments has been
conducted with the Lorenz-63 model for each fixed scaling parameter $c
\in [0.4,1.0]$ and a range of ensemble inflation factors. 
We only display the optimal results. The overall smallest RMS error is achieved
for $c = 0.6$ with a value of 4.1114.
The corresponding optimal RMS error for an ensemble square root filter is
$4.4813$ and $4.6596$ for the RHF, respectively. The EnKF with
perturbed observations leads to a RMS error of $4.1775$ which is only
slightly worse than the performance of the EGMF with $c=1$.}
\label{figN4.1}
\end{figure}

The three variable model 
\begin{equation}
\dot{x} = 10(y-x), \quad \dot{y} = x(28-z) - y, \quad \dot{z} = xy - \frac{8}{3} z
\end{equation}
of \cite{sr:lorenz63} is used as a final test for the EGMF method. 
Only the $x$ variable is observed every 0.20 time units with an observational error drawn from
a normal distribution with mean zero and variance eight. The model time-step is $\Delta t = 0.01$.
A total of 101000 assimilation steps is performed for each
experiment with only the last 100000 steps being used for the computation of RMS errors. 
We have implemented an ensemble square root filter, a RHF, and an EGMF
using formulation (\ref{FGM1})-(\ref{FGM2}) with ${\bf B} = c{\bf P}$, ${\bf P}$
the empirical covariance matrix of the ensemble. The parameter $c$ is
chosen from the interval $c \in [0.4,1.0]$. The number of ensemble
members is set to $M=25$ and no covariance localization is applied.  
The internal assimilation step-size is $\Delta s = 1/4$ and the 
$l_\infty$-norm of $u_{\rm B}({\bf
  x},s)$ is limited to a value of $u_{\rm cut} = 0.125/\Delta s$.
We have computed the RMS errors for
ensemble inflation factors  between 1.0 and 1.3 and only report the optimal results in
Fig.~\ref{figN4.1} as a function of $c$ for the EGMF. The overall smallest RMS error is achieved
for $c = 0.6$ with a value of 4.1114. The associated
RMS errors for the ensemble square root filter are $4.4813$ and $4. 6596$ for the RHF,
respectively. An increase in the number of ensemble members to $M=100$
leads to a reduction in the RMS error for the RHF to
4.3276 while the ensemble square root filter yields its optimal performance for $M = 50$ with a RMS
error of 4.3785. Both values are significantly 
larger than the optimal RMS error for the EGMF with $M=25$. A better
performance is observed for the EnKF with perturbed observations and
$M=25$ for which we obtain 4.1775 as the smallest RMS error which, in
fact, is close to the performance of the EGMF with $c=1$.

\section{Summary}

We have extended the popular EnKF to statistical models provided by
Gaussian mixtures.
The EGMF is derived using a continuous reformulation of the Bayesian analysis step and
consists of a combination of EnKF steps for each mixture component and
an exchange term. The exchange term is determined for each measurement
by a scalar elliptic PDE,  which can be solved analytically.
We have demonstrated by means of two numerical examples that the EGMF performs well when 
bimodal PDFs arise naturally from the model dynamics. The EGMF
provides a valuable and easy to implement alternative to sequential
Monte Carlo methods and other nonlinear filter algorithms. 
In this paper, we have used the standard EM algorithm
to assign Gaussian mixture models to ensemble predictions. More
refined methods such as those discussed by \cite{sr:fraley07} will be
considered in future work in order to provide a robust and accurate clustering 
of ensemble predictions.  Alternatively, one can implement the EGMF
with a Gaussian kernel density estimator. In this case, the empirical
covariance matrix of the ensemble can be used as a base for kernel 
bandwidth selection \citep{sr:wand}. With this choice, the EGMF becomes
closely related to the RHF of \cite{sr:anderson10}. The feasibility of our
approach has been demonstrated for the Lorenz-63 model. Further work
is required to assess the merits of Gaussian kernel density estimators
in comparison to EnKF and RHF implementations for high dimensional systems.
Encouraging results have also been reported by \cite{sr:stordal11} for
their related adaptive Gaussian mixture filter applied to the Lorenz-96 model 
\citep{sr:lorenz96,sr:lorenz98}. 


\bibliographystyle{plainnat}
\bibliography{survey}

\end{document}